\documentclass[a4paper, 12pt]{article}
\usepackage{bm}
\usepackage{amsfonts}
\usepackage{amsmath}
\usepackage{amssymb}
\usepackage{amsthm}
\usepackage{authblk}
\usepackage{color}
\usepackage{appendix}
\usepackage{listings}
\usepackage{xcolor}
\usepackage{textcomp}
\usepackage{slashbox}
\usepackage{url}

\usepackage{cite}
\usepackage[margin=2.0cm]{geometry}

\usepackage[utf8]{inputenc}
\usepackage[T1]{fontenc}

\usepackage{graphicx}
\usepackage[font=small,labelfont=bf]{caption}
\usepackage{subcaption}
\usepackage{float}

\usepackage{chngcntr}
\counterwithin{equation}{section}

\usepackage{array}
\allowdisplaybreaks

\begin{document}
	
	\title{Sum of Three Cubes via Optimisation}
	\author[1]{Boian Lazov\thanks{blazov\_fte@uacg.bg}}
	\author[2]{Tsvetan Vetsov\thanks{vetsov@phys.uni-sofia.bg}}
	\affil[1]{Department of Mathematics, University of Architecture, Civil Engineering and Geodesy, 1164 Sofia, Bulgaria}
	\affil[2]{Department of Theoretical Physics, Faculty of Physics, Sofia University, Sofia 1164, Bulgaria}
	
	\maketitle
	
	\begin{abstract}
		\noindent By first solving the equation $x^3+y^3+z^3=k$ with fixed $k$ for $z$ and then considering the distance to the nearest integer function of the result, we turn the sum of three cubes problem into an optimisation one. We then apply three stochastic optimisation algorithms to this function in the case with $k=2$, where there are many known solutions. The goal is to test the effectiveness of the method in searching for integer solutions. The algorithms are a modification of particle swarm optimisation and two implementations of simulated annealing. We want to compare their effectiveness as measured by the running times of the algorithms. To this end, we model the time data by assuming two underlying probability distributions -- exponential and log-normal, and calculate some numerical characteristics for them. Finally, we evaluate the statistical distinguishability of our models with respect to the geodesic distance in the manifold with the corresponding Fisher information metric.
	\end{abstract}
	
	
	\section{Introduction}
	The sum of three cubes problem can be stated in the following way: Let $k$ be a positive integer. Is there a solution to the equation
	\begin{align}
	x^3+y^3+z^3=k,\label{deq}
	\end{align}
	such that $(x,y,z)\in \mathbb{Z}^3$?

	As is well known, for $k\equiv 4\ (\mathrm{mod}\ 9)$ such a solution does not exist \cite{davenport39}. For $k\not\equiv4\ (\mathrm{mod}\ 9)$ however, it has been conjectured by Heath-Brown that there are infinitely many solutions \cite{heath-brown92}. A direct search for solutions is one way to support this conjecture. Until recently there were only two numbers below $100$ for which a representation as a sum of three integer cubes had not been found. Those were $33$ and $42$.
	
	Then the sum of three cubes problem gained an unusual amount of fame in the last year after a Numberphile video inspired a solution to the case with $k=33$ \cite{booker19}. What followed were a solution for $k=42$ and a new (third) solution for $k=3$.
	
	The main problem when directly searching for solutions is the time it takes for a brute force approach. However, there are ways to reduce this time and the latest method by Booker reached a time complexity of $O\left(B^{1+\varepsilon}\right)$, with $\mathrm{min}\left\{ |x|,|y|,|z| \right\}\le B$ \cite{booker19}, i.e. an almost linear search.
	
	When thinking about a way to improve this result, a natural step seems to be to gamble a bit and rely on the odds of probability. In other words we may try to guess the solution via some random search heuristic. One way this can be achieved is by turning the sum of three cubes into an optimisation problem, for which there are a wealth of such heuristics. As it turns out there have been some attempts to solve diophantine equations using particle swarm optimisation (PSO) algorithms \cite{abraham10}. And while we use PSO in some form, our approach is different and permits the use of many stochastic optimisation methods.
	
	This paper is organised as follows. In section \ref{meth} we define the function to be optimised. Then, since we've said that our approach was motivated by the desire to use a random search heuristic, in the next sections we apply to our function two such heuristics that look promising, namely PSO and simulated annealing (SA). Their performance is then evaluated experimentally by measuring the time it takes the respective algorithm to obtain a solution to (\ref{deq}) for $k=2$. This choice of $k$ is justified by its high density of solutions, which makes testing much easier. In particular, in section \ref{dfo} we introduce our version of a PSO algorithm, while in section \ref{sa} we use two versions of SA -- one with restarts and one without. In section \ref{stat} we conduct a thorough statistical analysis of the performance of our algorithms and of their similarities. We conclude this work with section \ref{concl}, where we briefly comment on our results and show some possible directions for future research.
	
	\section{Our approach}\label{meth}
	
	We first start by defining the function, which we will be optimizing. Solving equation (\ref{deq}) for $z$ we trivially get
	\begin{align}
	z=\left( k-x^3-y^3 \right)^{\frac{1}{3}}.\label{zeq}
	\end{align}
	We can then define a function $f_k(x,y)$ as
	\begin{align}\label{fdef}
	f_k(x,y):=\left\| z \right\|
	= \left\| \left( k-x^3-y^3 \right)^{\frac{1}{3}} \right\|,
	\end{align}
	where $\left\Vert z \right\Vert$ denotes the distance to the nearest integer from $z$.\footnote{For a summary of the main properties of the distance to the nearest integer function see \url{https://www.researchgate.net/publication/308023356_Note_on_the_Distance_to_the_Nearest_Integer}}
	
	Now let us fix $k=k_0$. If $(x_0,y_0,z_0)$ is an integer solution to Eq. (\ref{deq}), we have
	\begin{align}
	f_{k_0}(x_0,y_0)=0\label{min}
	\end{align}
	and this is a global minimum. Conversely, if $f_{k_0}(x,y)$ has a local minimum at $(x_0,y_0)\in \mathbb{Z}^2$, such that (\ref{min}) is satisfied, it has a global minimum there and this gives the solution $(x_0,y_0,z_0)$ to (\ref{deq}), where $z_0$ is evaluated from Eq. (\ref{zeq}).
	
\begin{figure}[H]
	\centering
	\begin{subfigure}[b]{0.51\textwidth}
		\includegraphics[width=\textwidth]{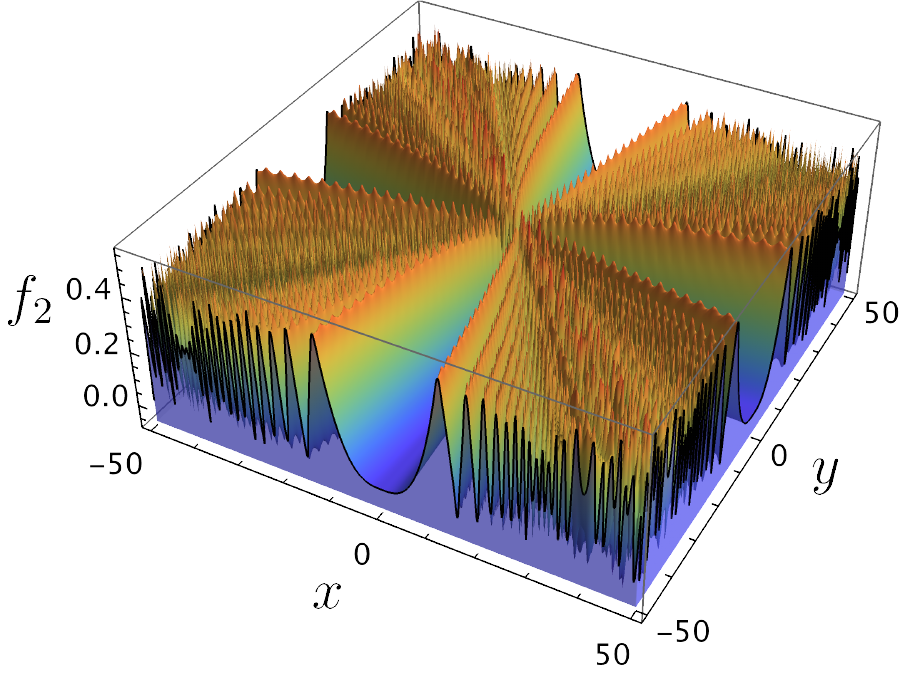}
		\caption{Many local optima of the function $f_2(x,y)$.}
		\label{fig:fPlot1}
	\end{subfigure}
	\hspace{0.8 mm}
	\begin{subfigure}[b]{0.46\textwidth}
		\includegraphics[width=\textwidth]{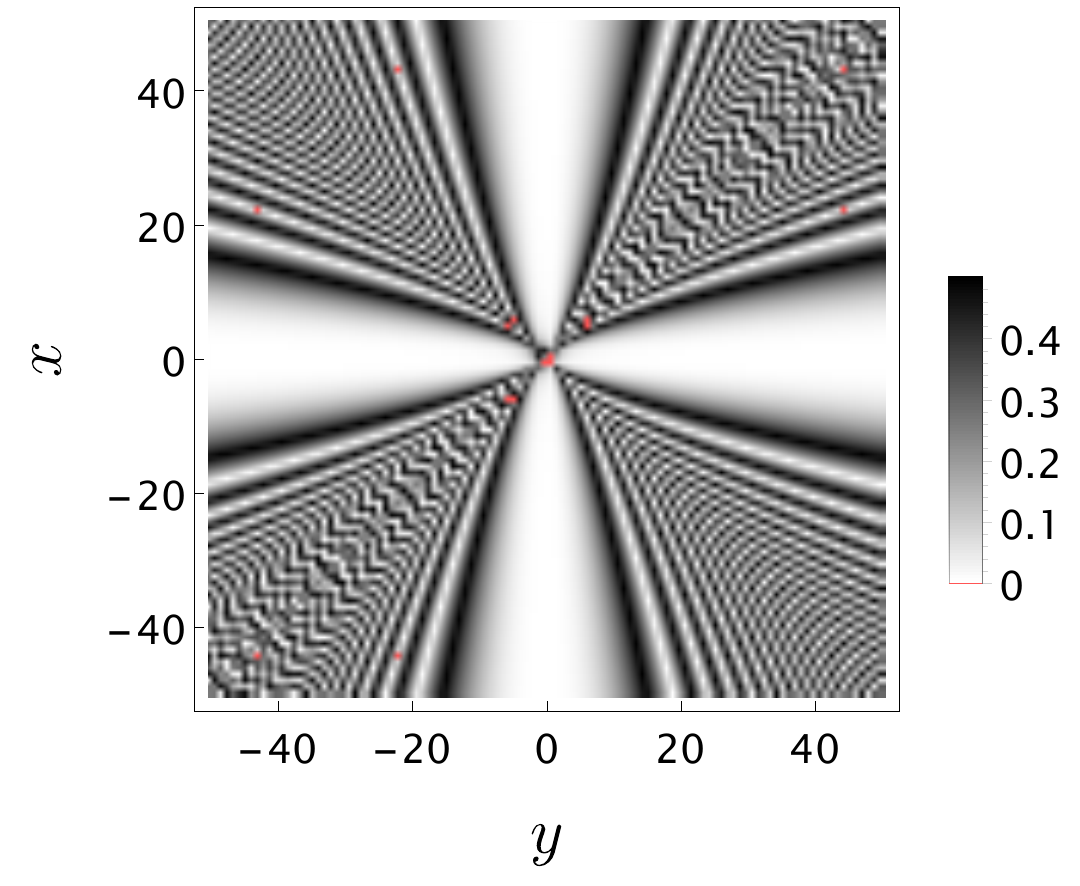}
		\caption{A discrete matrix plot of $f_2(x,y)$.}
		\label{fig:fPlot2}
	\end{subfigure}
\caption{Sample of the function $f_2(x,y)$ in the range $x,y\in[-50,50]$. (a) depicts the many local optima of the function $f_2(x,y)$. (b) shows a discrete matrix plot of $f_2(x,y)$ with red pixels corresponding to the global minima, where $f_2(x,y)=0$.}\label{fig:fPlot}
\end{figure}

	Hence, the problem now is to find a global minimum of $f_k(x,y)$ with $(x,y)\in \mathbb{Z}^2$ and fixed $k$. As it turns out, this is not an easy problem since our function has no shortage of local optima. As an illustration we show a plot of the function $f_2(x,y)$ with $x,y\in[-50,50]$ in figure \ref{fig:fPlot}.

	\section{Swarm optimisation}\label{dfo}
	
	The first algorithm we use to minimise (\ref{fdef}) is a representative from the family of swarm algorithms. These algorithms try to find extrema of a function by emulating the behaviour of swarms of insects searching for resources. There are many algorithms based on this idea and they may differ significantly in performance, depending on the problem. For the latest development on the subject one may refer to \cite{cuevas:2020}.
	
	For our particular function we first tried using standard particle swarm optimisation (SPSO) \cite{clerc12}. However, it turned out to be too slow and, searching for a better alternative, we stumbled upon dispersive flies optimisation (DFO) developed by Al-Rifaie \cite{al-rifaie14}. It has the benefit of simplicity and we decided to use it as a base and see where we end up. As it turned out we made some significant changes in order to improve the performace for our particular function.
	
	As the main point we found the breaking of the swarm in DFO to be too limited for our needs, so we opted for simple randomisation of all flies after a certain criterion has been met. This, however, meant that the swarm best would change after every dispersion and may have lead to insufficient exploration around it. To fix this we introduced a memory of the best position of the swarm found so far. We then used this ``best swarm best'' in the position update formula in addition to the neighbours best and swarm best positions.
	
	This improved the performance of our algorithm when compared with unmodified DFO. However, our function has a lot of local minima and the algorithm tended to get stuck in them often. Furthermore, the closer the best swarm best position's fitness function was to $0$, the harder it was to find a better position (including after a dispersion). One way to alleviate this somewhat was to introduce simple restarts after a fixed number of iterations without an improvement in the best swarm best position.
	
		This required including a new parameter in the algorithm however, namely the number of iterations. To avoid it we decided to instead probabilistically change the best swarm best to the current swarm best in the iteration even if the new position was worse. The probability depends on the difference between the fitness function values in both positions.
	
	\subsection{Description of the algorithm}

	Our PSO algorithm can in principle be used for any discrete optimisation problem and in this subsection we will give a general description, which doesn't refer to specifics such as the dimension of the search space or the explicit form of the fitness function.
	
	Like all PSO algorithms it first starts with the initialisation of the swarm (with $s$ particles). This is done by choosing a random position $\vec{x}_i$ for each particle $i\in \{0,...,s-1\}$ inside the (discrete) search space and calculating it's fitness function $f\!f(\vec{x}_i)$. Then the swarm best position $\vec{sb}\in\{ \vec{x}_0,...,\vec{x}_{s-1} \}$ is determined, such that
	\begin{align}
	f\!f(\vec{sb})=\min_{i\in\{0,...,s-1\}}\{ f\!f(\vec{x}_i) \}.
	\end{align}
	
	As mentioned, the algorithm uses a memory of the best position obtained so far, which we call best swarm best $\vec{bsb}$. At initialisation this is set equal to the swarm best, i.e. $\vec{bsb}\leftarrow \vec{sb}$.\\
	\indent As usual, the particles communicate with their neighbours. We use the standard ring topology for the set of neighbours $n_i$ of particle $i$ \cite{clerc12}, i.e.
	\begin{align}
	n_i=\left\{( i-1)\  \mathrm{mod}\  s,\ i,\ (i+1)\ \mathrm{mod}\ s \right\}.
	\end{align}
	Then, knowing the neighbours, the best position among them is found for each particle. We call this the neighbours best $\vec{nb}_i$ and it satisfies the following 
	\begin{align}
	f\!f(\vec{nb}_i) = \min_{k\in n_i}\{f\!f(\vec{x}_k)\},\ \ \ \vec{nb}_i\in \{ \vec{x}_k \}_{k\in n_i}.
	\end{align}
	
	Next follow the iterations. We can limit the number of iterations to get some approximate solution or wait for some condition to be satisfied. An iteration consists of a position update, confinement and a $\vec{sb}$, $\vec{bsb}$ and $\vec{nb}_i$ update.
	
	First is the position update. Since, as we've said, the swarm is periodically dispersed, the position update formula depends on a simple dispersion condition. We first define a dispersion parameter $dp=s/5$. Then the number of particles that have reached the best swarm best position is counted and, if they are no less than $dp$, the positions of all particles are randomised across the search space.
	
	If the dispersion condition has not been met, the positions of the particles are updated by the formula
	\begin{align}
	x_{i,d}\leftarrow \mathrm{round}\left( nb_{i,d} + \frac{r}{2}\left(sb_d+bsb_d-2 x_{i,d}\right) \right),\ \ \ d=1,...,D.
	\end{align}
	where $r$ is drawn from a uniform distribution on the interval $(0,1)$ and $D$ is the dimension of the search space.
	
	After the position update there may be some particles outside the search space. We want them confined inside however, so random positions are chosen for such particles. While this is a simple way to implement confinement, it helps with the exploration behaviour that some problems so desperately need.
	
	Finally, $\vec{sb}$, $\vec{bsb}$, and $\vec{nb}_i$ need to be updated. After $\vec{sb}$ is determined as in the initialisation, it is used to update $\vec{bsb}$. If $f\!f(\vec{sb})<f\!f(\vec{bsb})$,
	\begin{align}
	\vec{bsb}\leftarrow \vec{sb},
	\end{align}
	as expected. If, however, $f\!f(\vec{sb})>f\!f(\vec{bsb})$, the worse position $\vec{sb}$ is accepted as the new best swarm best with probability
	\begin{align}
	p=1-\frac{f\!f(\vec{sb})-f\!f(\vec{bsb})}{0.5}.
	\end{align}
	\indent As can be seen, the choice of $\vec{bsb}$ borrows it's idea from SA algorithms. In practice the change to a worse $\vec{bsb}$ happens after the particles have dispersed, because this is when $f\!f(\vec{sb})$ can be less than $f\!f(\vec{bsb})$.
	
	Next, $\vec{nb}_i$ is determined as in the initialisation. Naturally, after each $\vec{bsb}$ update it needs to be checked whether some suitable condition has been satisfied so the algorithm can be exited.
	
	\subsection{Computational results}\label{results}
	
	In order to experiment with our PSO algorithm we need to fix some parameters for our particular problem. First of all, we want to search for solutions to the diophantine equation (\ref{deq}) with $k=2$, i.e. global minima of $f_2(x,y)$. This means that our search space has $2$ dimensions, i.e. $D=2$. So we denote the positions in the search space with $(x,y)$ and we fix our fitness function as
	\begin{align}
	f\!f(x,y)=f_2(x,y)
	=\left\| \left( 2-x^3-y^3 \right)^{\frac{1}{3}} \right\|.
	\end{align}
	
	We also need to determine the exit criterion. As we are searching for integer solutions and thus an approximate one is not good enough, we first choose some threshold $thr$ and the algorithm looks for a pair $(x_0,y_0)$, such that $f\!f(x_0,y_0)\le thr$. Then the candidate solution $\left(x_0,y_0,\mathrm{round}\left(z_0\right)\right)$, where $z_0$ is calculated from (\ref{zeq}), is plugged into (\ref{deq}) and the algorithm is exited, if the equation is satisfied.
	
	Additionally, after some experimenting, that is in no way conclusive, we found that a particle swarm size of $s=50$ works best in our case.\\
\indent With the above fixed we tested the time performance of our algorithm for different ranges of $x$ and $y$ by recording the time it needs to find a solution. The full code which we used for testing is included in appendix \ref{psocode}. We wrote the algorithm in C.\\
\indent We chose to scan three different ranges of values for $x$ and $y$, namely
\begin{align}\label{eqRanges}
\nonumber &R_3=\{(x,y):\ x,y\in\mathbb{Z},\ -10^3\leq x\leq0\leq y\leq 10^3\},\quad N=10^4,\\
  &R_4=\{(x,y):\ x,y\in\mathbb{Z},\ -10^4\leq x\leq0\leq y\leq 10^4\},\quad N=10^4,\\
 \nonumber   &R_5=\{(x,y):\ x,y\in\mathbb{Z},\ -10^5\leq x\leq0\leq y\leq 10^5\},\quad N=10^3.
\end{align}
Here we denote with $N$ the respective number of runs completed by the algorithm in the corresponding range. Hence, $N$ is also the sample size of the time data for a given range.

\begin{figure}[H]
		\centering
		\begin{subfigure}[b]{0.47\textwidth}
			\includegraphics[width=\textwidth]{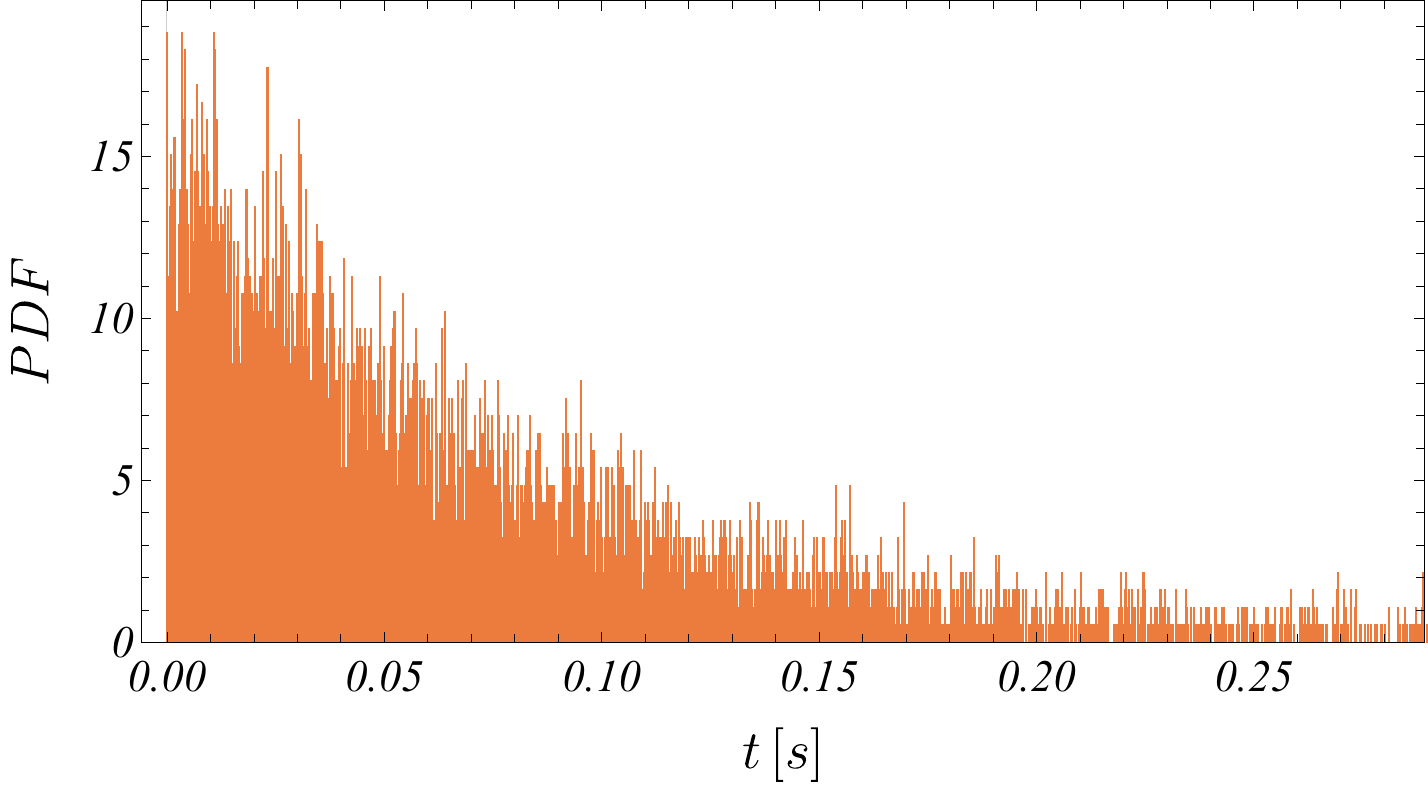}
			\caption{PSO in the range $R_3$}
			\label{fig:DFO-Hist1}
		\end{subfigure}
		\hspace{0.8 mm}
		\begin{subfigure}[b]{0.47\textwidth}
			\includegraphics[width=\textwidth]{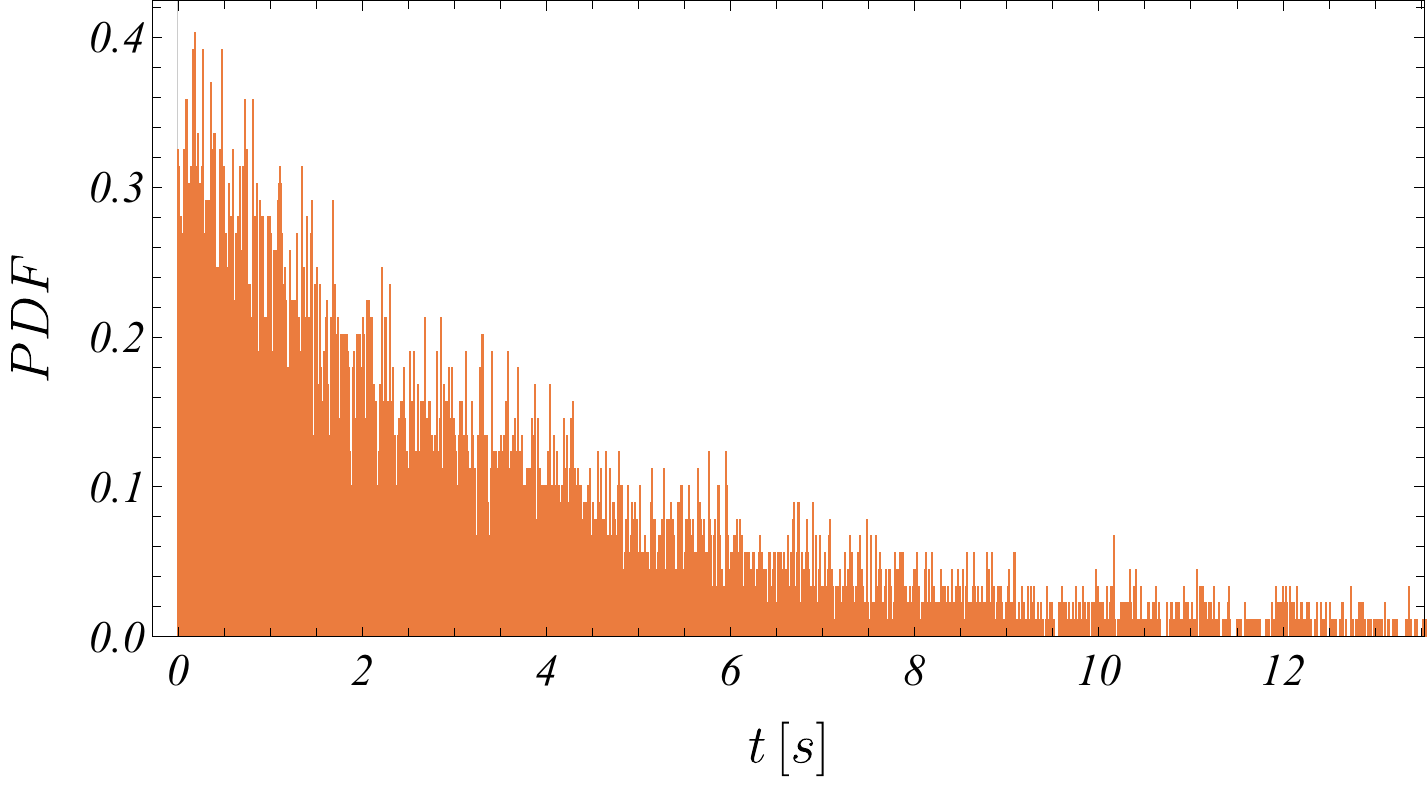}
			\caption{PSO in the range $R_4$}
			\label{fig:DFO-Hist2}
		\end{subfigure}
		\hspace{0.8 mm}
		\begin{subfigure}[b]{0.47\textwidth}
			\includegraphics[width=\textwidth]{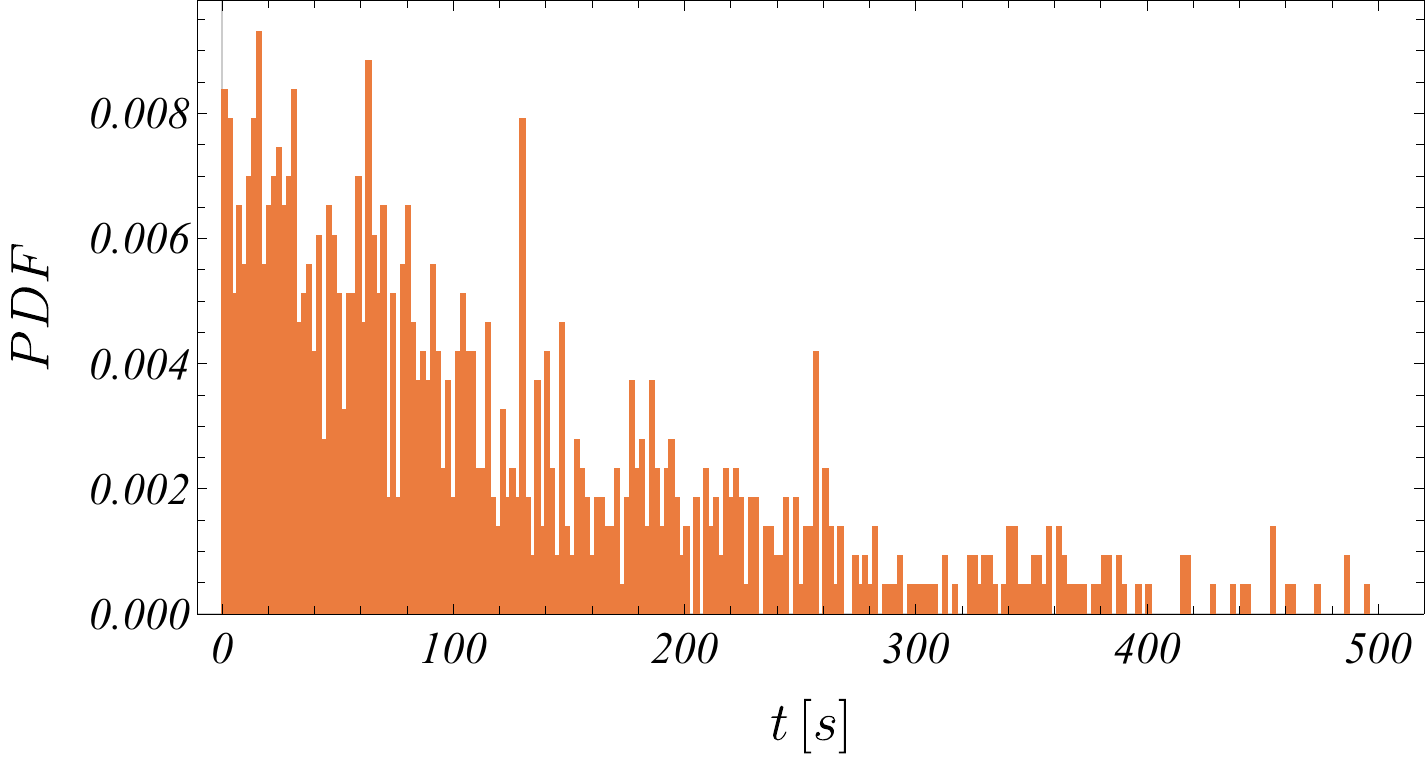}
			\caption{PSO in the range $R_5$}
			\label{fig:DFO-Hist3}
		\end{subfigure}
		\caption{Time performance of PSO. Every histogram gives the distribution of the individual times the PSO algorithm took to find a solution in the: (a) range $R_3$ with sample size $N=10^4$; (b)  range $R_4$ with sample size $N=10^4$; (c) range $R_5$ with sample size $N=10^3$.}\label{fig:DFOH}
	\end{figure}

The results are arranged in histograms based on the data set $\{t_i\}_{i=1}^N$ of running times $t_i$. The data is gathered into bins with equal widths
\begin{equation}
\Delta t=\frac{{\rm{max}} \{t_i\}- {\rm{min}} \{t_i\}}{\ell},
\end{equation}
where ${\rm{max}} \{t_i\}$ is the longest time for finding a solution, ${\rm{min}} \{t_i\}$ is the shortest one, and $\ell<N$ is an arbitrary partition. The data is also normalised to show the corresponding probability density function (PDF) (figure \ref{fig:DFOH}). In this case, the data is discrete and the PDF function states that the probability of $t_i$ falling within an interval of width $\Delta t$ is the density in that range times $\Delta t$.

	\section{Simulated annealing}\label{sa}
	
	SA is a search heuristic based on the physical process of annealing \cite{bertsimas93}. It has the advantage of being effective despite its ease of implementation. SA has been extensively studied and has many variations, e.g. in the choice of cooling schedule or neighbours. Those variations of course affect the algorithm's performance.
	
	One particularly interesting modification to SA is the inclusion of restarts. It has been shown that under some conditions restarting the algorithm according to certain criteria results in improved times for finding a desired extremum \cite{mendivil01}. More precisely a restarting SA (rSA) algorithm with a local generation matrix and cooling schedule $\mathrm{temp}(m)\sim\frac{1}{m}$ has probabilities that the extremum has not been reached by time $m$ which converge to zero at least geometrically fast in $m$.
	
	Here we implement two versions of SA -- one without restarts and one with restarts. As will be seen, we are using a logarithmic cooling schedule and thus our implementation of the algorithm fails to satisfy all the assumptions of theorem 4.1 of \cite{mendivil01}. Nevertheless, those are not necessary conditions and it turns out that restarting significantly improves the running time in our particular case.

\subsection{Computational results}

	As in section \ref{results}, we again want to test the performance of the algorithm for $f_2(x,y)$, so this is our energy function in the context of SA. A state here is just a point $(x,y)\in\mathbb{Z}^2$. The exit criterion is also the same as in section \ref{results}. We use the following cooling schedule:
\begin{align}
\mathrm{temp}(m)=\frac{1}{\ln m+0.01},
\end{align}
where $m$ as usual is the current iteration number.

The neighbourhood of a state $(x,y)$ is
\begin{align}
n(x,y)=\left\{ (x+a,y+b) : a,b\in\mathbb{Z},\ a,b\in[-10,10] \right\}
\end{align}
except for the states close to the border of the search space, where we remove the appropriate points so as not to end up outside. The generation matrix allows transitions from $(x,y)$ to all points in $n(x,y)$ with equal probability except to $(x,y)$ itself.

\begin{figure}[H]
	\centering
	\begin{subfigure}[b]{0.47\textwidth}
		\includegraphics[width=\textwidth]{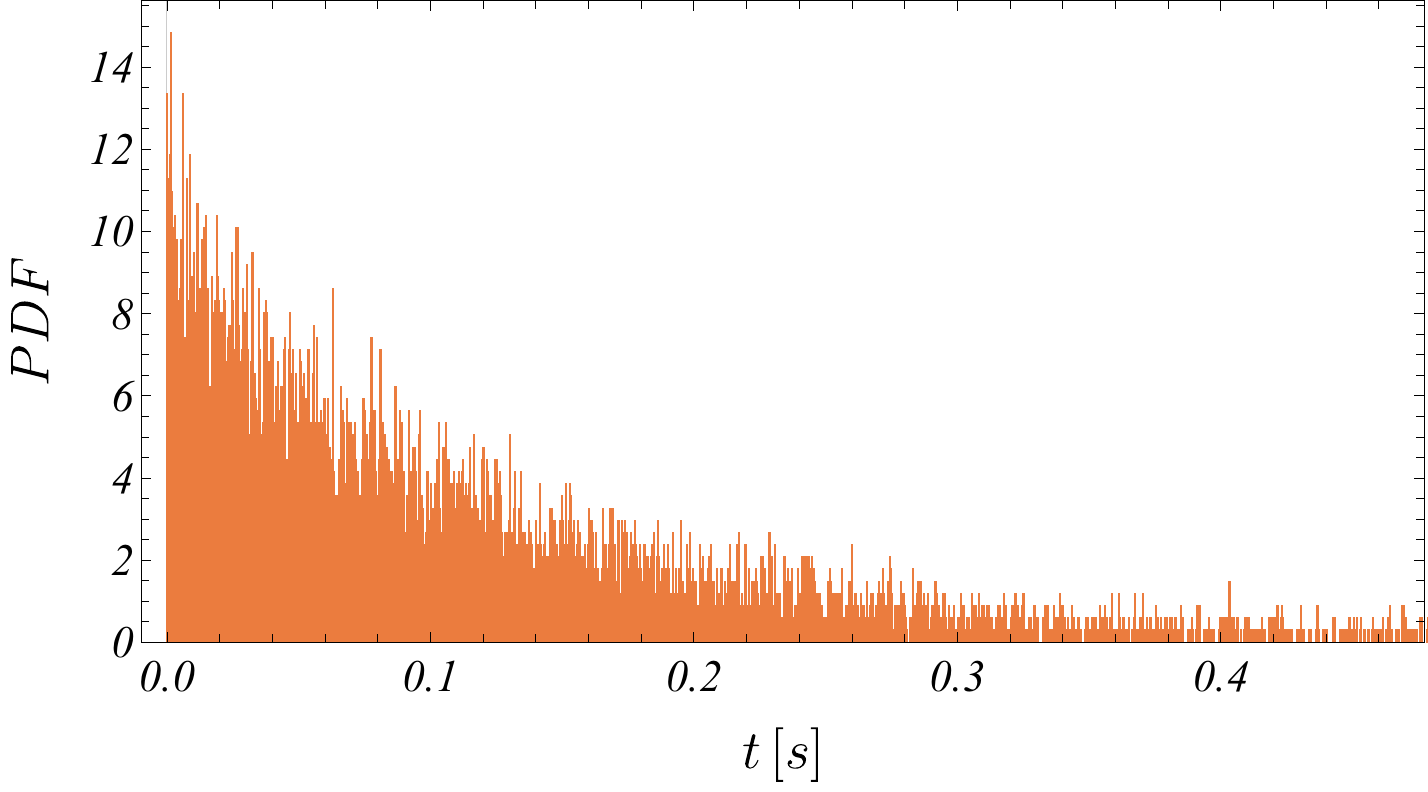}
		\caption{SA in the range $R_3$}
		\label{fig:SANORES-Hist1}
	\end{subfigure}
	\hspace{0.8 mm}
	\begin{subfigure}[b]{0.47\textwidth}
		\includegraphics[width=\textwidth]{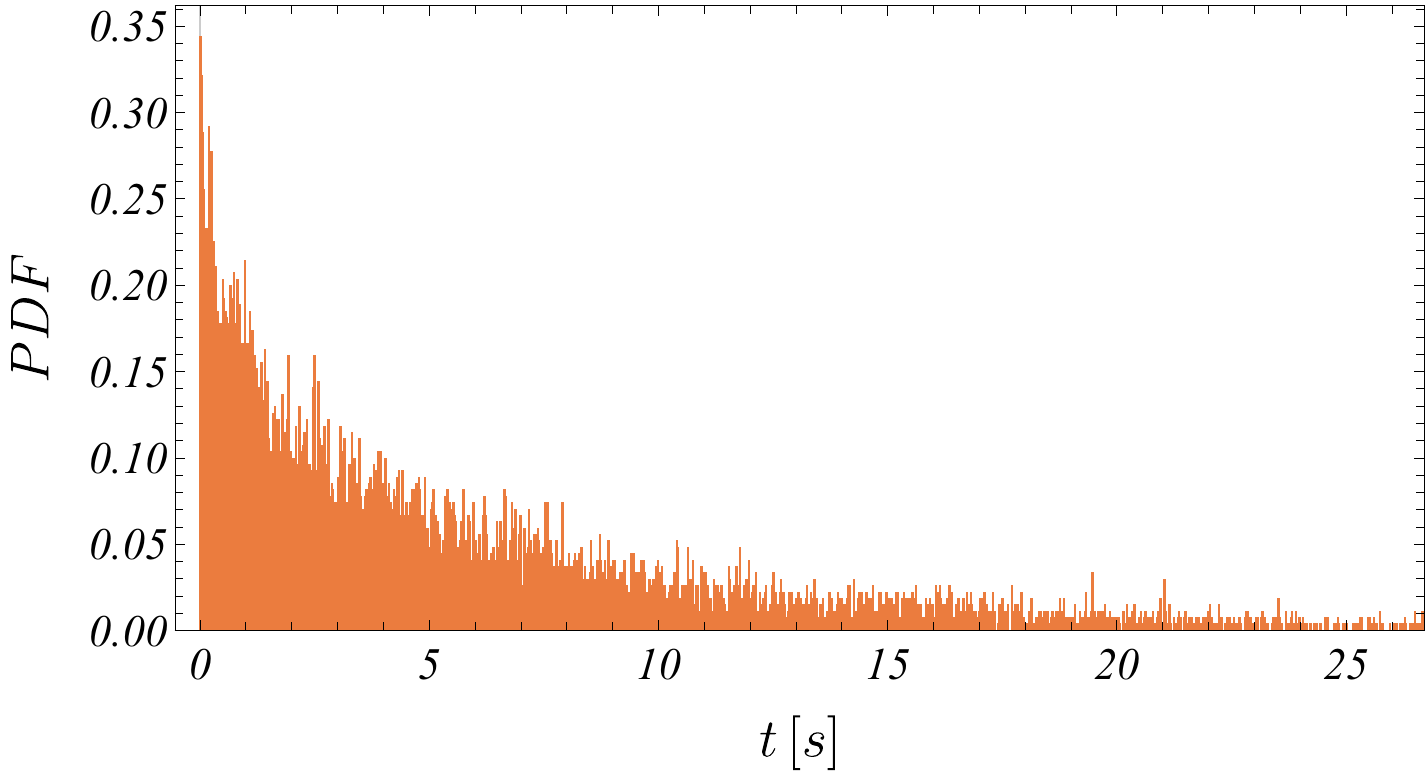}
		\caption{SA in the range $R_4$}
		\label{fig:SANORES-Hist2}
	\end{subfigure}
	\hspace{0.8 mm}
	\begin{subfigure}[b]{0.47\textwidth}
		\includegraphics[width=\textwidth]{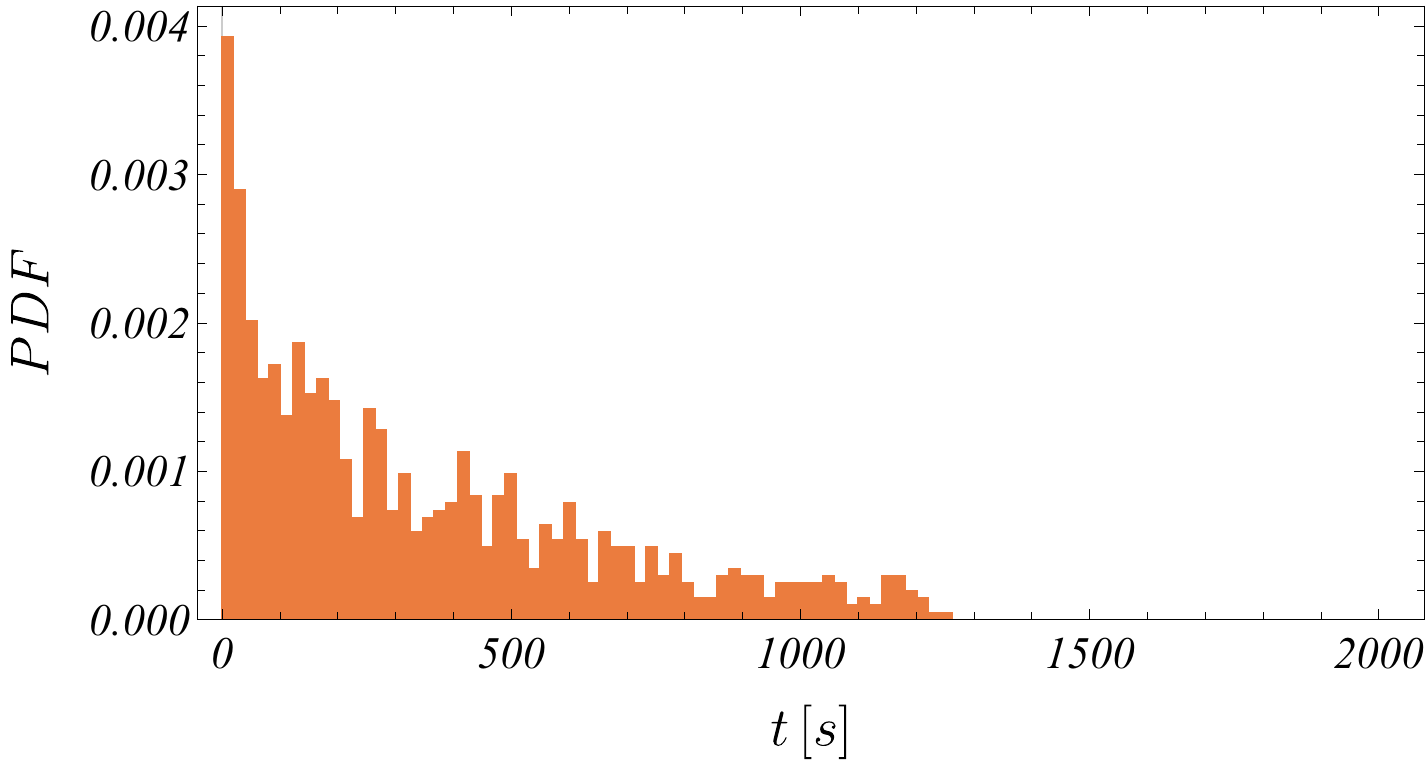}
		\caption{SA in the range $R_5$}
		\label{fig:SANORES-Hist3}
	\end{subfigure}
	\caption{Time performance of SA. Every histogram gives the distribution of the individual times the SA algorithm took to find a solution in the: (a) range $R_3$ with sample size $N=10^4$; (b)  range $R_4$ with sample size $N=10^4$; (c)  range $R_5$ with sample size $N=10^3$.}\label{fig:SANORESH}
\end{figure}

\begin{figure}[H]
	\centering
	\begin{subfigure}[b]{0.47\textwidth}
		\includegraphics[width=\textwidth]{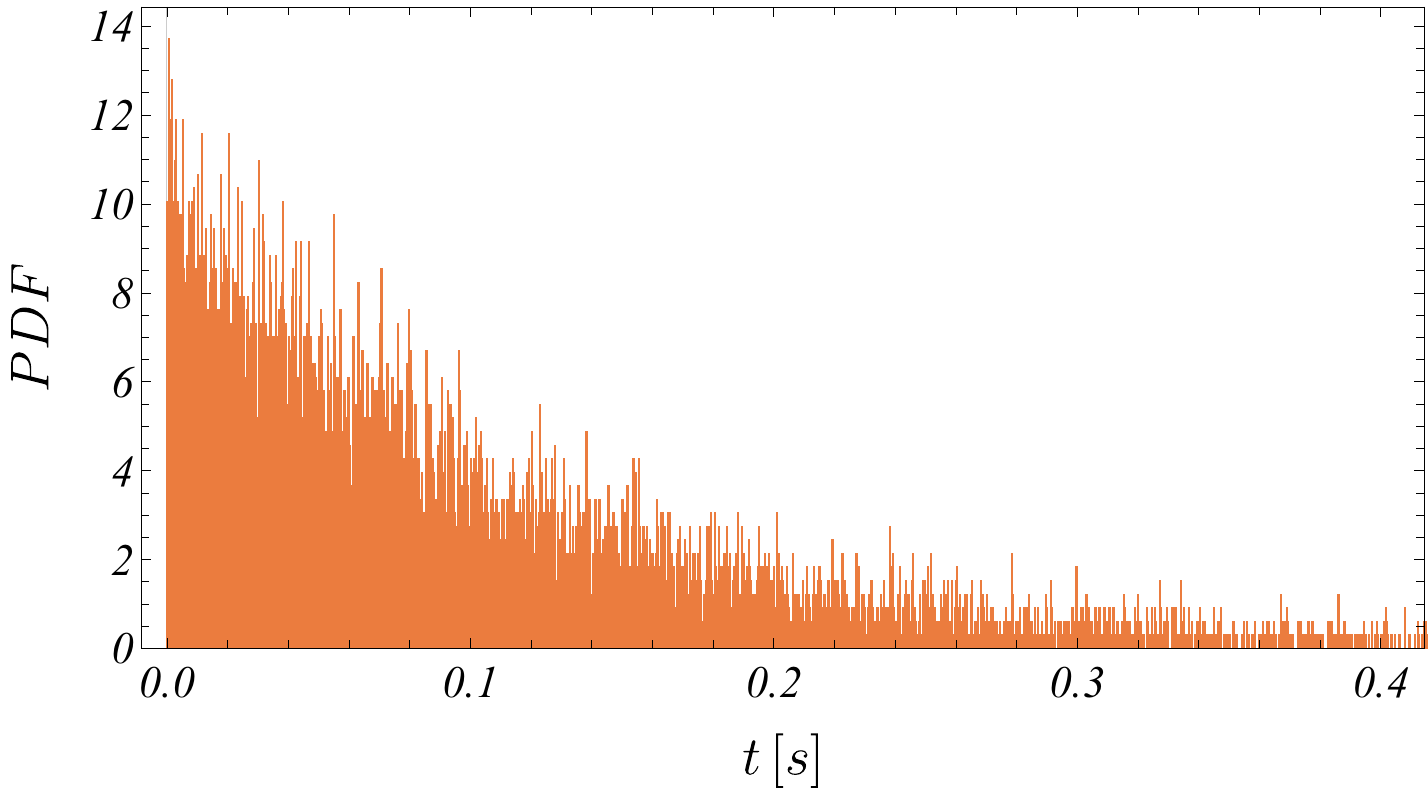}
		\caption{rSA in the range $R_3$}
		\label{fig:SARES-Hist1}
	\end{subfigure}
	\hspace{0.8 mm}
	\begin{subfigure}[b]{0.47\textwidth}
		\includegraphics[width=\textwidth]{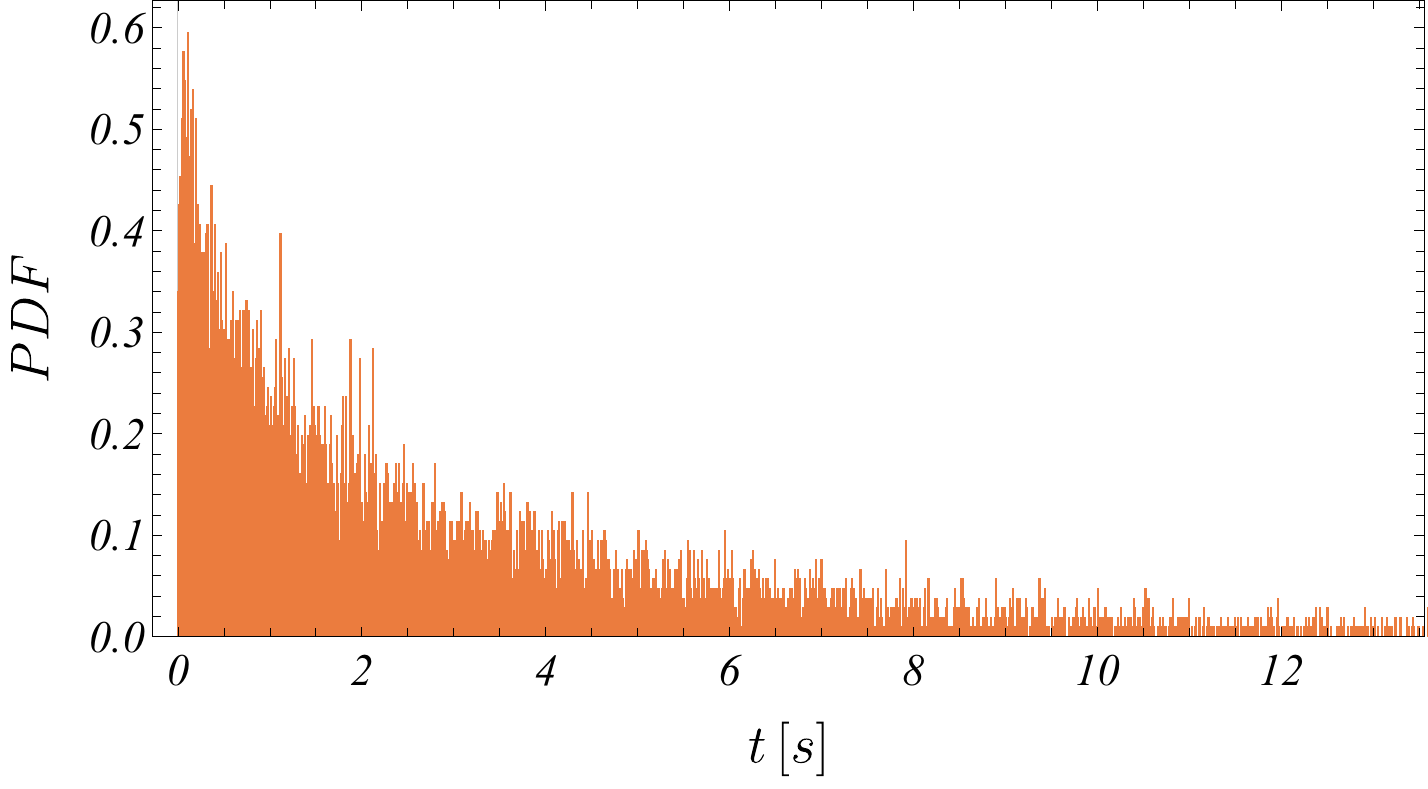}
		\caption{rSA in the range $R_4$}
		\label{fig:SARES-Hist2}
	\end{subfigure}
	\hspace{0.8 mm}
	\begin{subfigure}[b]{0.47\textwidth}
		\includegraphics[width=\textwidth]{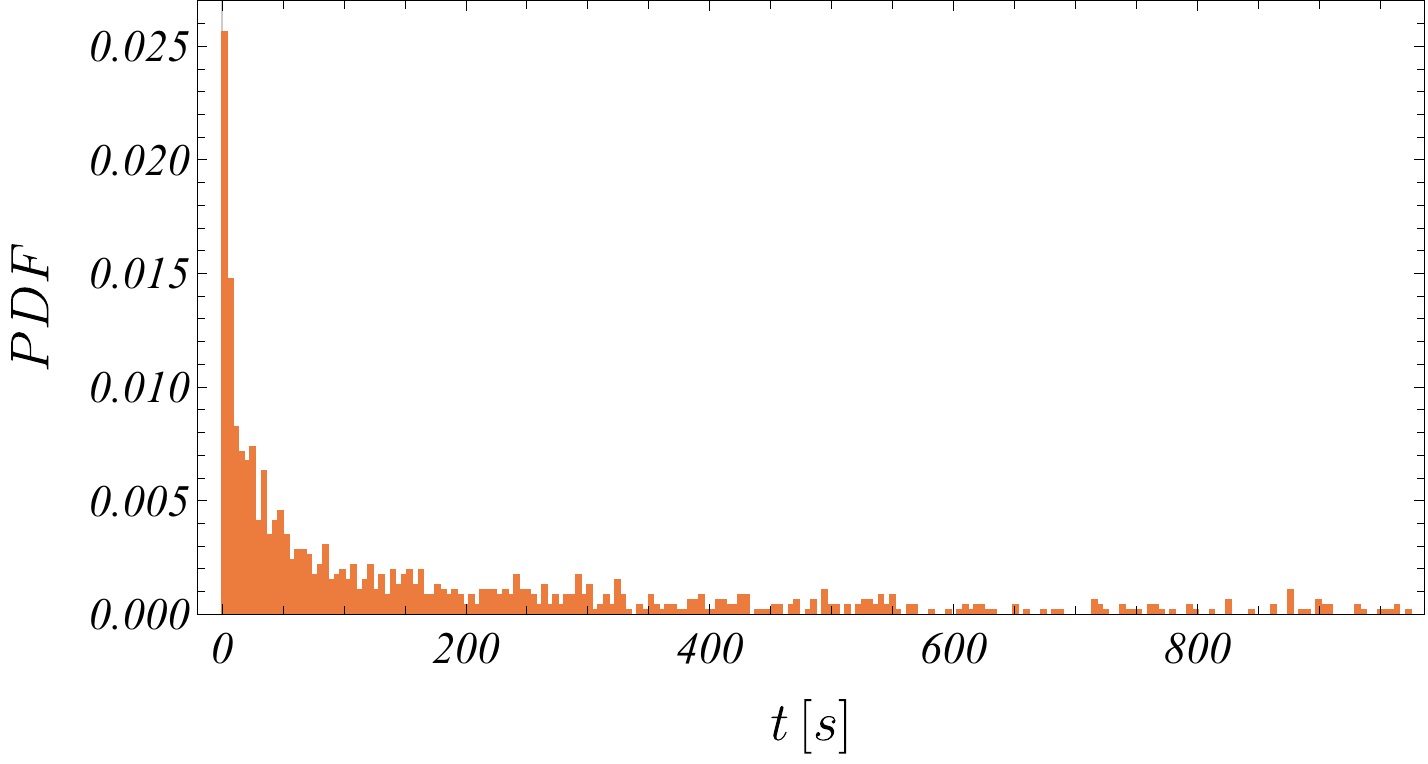}
		\caption{rSA in the range $R_5$}
		\label{fig:SARES-Hist3}
	\end{subfigure}
	\caption{Time performance of rSA. Every histogram gives the distribution of the individual times the rSA algorithm took to find a solution in the: (a) range $R_3$ with sample size $N=10^4$; (b)  range $R_4$ with sample size $N=10^4$; (c) range $R_5$ with sample size $N=10^3$.}\label{fig:SARESH}
\end{figure}

For the rSA algorithm, we use the criterion suggested in \cite{mendivil01}, namely restarting after $rtm$ consecutive states have the same energy. After some experimenting, we decided to use $rtm=30$.\\
\indent Again, we wrote the algorithms in C and include the full codes in apprendices \ref{sacode} and \ref{rsacode}. Figures \ref{fig:SANORESH} and \ref{fig:SARESH} show the PDF histograms for the two versions of SA in different ranges.

\section{Data models and statistical analysis}\label{stat}

In this section we conduct a standard statistical analysis by modeling the time data, produced by the given algorithms, with carefully chosen continuous probability distributions. Our goal is to evaluate the relative performance of the three algorithms and furthermore to estimate their similarities. We chose to describe the accumulated time data by two statistical models, namely  a simple one-parameter exponential model $f(t;\lambda)$, and a two-parameter log-normal distribution $f(t;\alpha,\beta)$. Here the stochastic variable is the individual time $t$ a given method takes to find an integer solution to $x^3+y^3+z^3=2$. 

\subsection{Exponential distribution}

The exponential model is a simple one-parameter probability distribution, where one assumes that the underlying statistics models a Poisson process. The PDF of an exponential distribution is given by
\begin{equation}\label{eqExpD}
f(t;\lambda ) = \left\{ \begin{array}{l}
\lambda {e^{ - \lambda t}},\,\,t \ge 0\\
0,\,\,t < 0
\end{array} \right.,
\end{equation}
where $\lambda>0$ is the rate parameter. 

The expected value, the variance and the median of an exponentially distributed random variable $t$ with rate parameter $\lambda$ are well known, namely
\begin{equation}
\bar t={\mathop{\rm E}\nolimits} [t] = \frac{1}{\lambda },\quad {\mathop{\rm Var}\nolimits} [t] = \frac{1}{{{\lambda ^2}}},\quad {\mathop{\rm Med}\nolimits}[t] = \frac{{\ln 2}}{\lambda }.
\end{equation}

When a finite sample data is available the mean time $\bar{t}$ for finding a solution also coincides with the mean time from the sample data:
\begin{equation}\label{eqMeanTSample}
\bar{t}=\frac{1}{{N}}\sum\limits_{i=1}^N{t_i}.
\end{equation}

The $95\% $ confidence intervals for $\lambda$ and $t$ are given by  
\begin{equation}
{\lambda _{lower}} \le \lambda  \le {\lambda _{upper}},\quad \frac{1}{{{\lambda _{upper}}}} \le \bar t \le \frac{1}{{{\lambda _{lower}}}},
\end{equation}
where
\begin{equation}
{\lambda _{lower}} = \lambda \left( {1 - \frac{{1.96}}{{\sqrt N }}} \right),\quad {\lambda _{upper}} = \lambda \left( {1 + \frac{{1.96}}{{\sqrt N }}} \right).
\end{equation}

\subsection{Log-normal distribution}

The two-parameter log-normal distribution $f(t;\alpha,\beta)$, $\alpha\in(-\infty,\infty)$, $\beta>0$, is a continuous probability distribution of a positive random variable $t>0$, whose logarithm is normally distributed. There are many different parameterisations of the log-normal distribution, but we prefer the following:
\begin{equation}\label{eqLogND}
f(t;\alpha ,\beta ) = \frac{1}{{\beta \sqrt {2\pi } t}}{{\rm{e}}^{ - \frac{{{{(\ln t - \alpha )}^2}}}{{2{\beta ^2}}}}},
\end{equation}
where the parameters of the distribution can be obtained directly from the sample data via
\begin{equation}
\alpha  = \frac{1}{N^{}}\sum\limits_{i = 1}^{N^{}} {\ln {t_i}},\quad \beta  = \frac{1}{{\sqrt {N^{}} }}\sqrt {\sum\limits_{i = 1}^{N^{}} {{{(\ln {t_i} - \alpha )}^2}} }.
\end{equation}

In this case, the mean time $\bar{t}$ for finding a solution and its standard deviation are given by
\begin{equation}
\bar t = {\mathop{\rm E}\nolimits} [t] = {e^{\alpha  + \frac{{{\beta ^2}}}{2}}}\quad {\mathop{\rm SD}\nolimits} [t]=\sqrt{{\mathop{\rm Var}\nolimits} [t]} = {e^{\alpha  + \frac{\beta ^2}{2}}}\sqrt{{e^{{\beta ^2}}} - 1},
\end{equation}
where $\bar t$ does not coincide with the mean sample time (\ref{eqMeanTSample}). Furthermore, the median and the mode yield
\begin{equation}
{\mathop{\rm Med}\nolimits} [t] = {e^\alpha },\quad {\mathop{\rm Mode}\nolimits} [t] = {e^{\alpha  - {\beta ^2}}} ,
\end{equation}
where the mode defines the point of global maximum of the probability density function. 

The standard scatter intervals for the log-normal distribution are written by 
\begin{eqnarray}
{t_{68\% }} \in [{e^{\alpha  - \beta }},{e^{\alpha  + \beta }}], \quad {t_{95\% }} \in [{e^{\alpha  - 2\beta }},{e^{\alpha  + 2\beta }}].
\end{eqnarray}
However, these estimates are not very informative for skew-symmetric distributions such as the log-normal one. In this case, we can extract an efficient $95\%$ confidence interval for the log-normal model based on the Cox proposal, namely \cite{Zhou}
\begin{equation}\label{eqCoxConfInt1}
{\bar t_{95\% }} \in {e^{\left[ {\alpha  + \frac{{{\beta ^2}}}{2} - 1.96\sqrt {\frac{{{\beta ^2}}}{{{N^{}}}} + \frac{{{\beta ^4}}}{{2({N^{}} - 1)}}} ,{\kern 1pt} \alpha  + \frac{{{\beta ^2}}}{2} + 1.96\sqrt {\frac{{{\beta ^2}}}{{{N^{}}}} + \frac{{{\beta ^4}}}{{2({N^{}} - 1)}}} } \right]}},
\end{equation}
where we can estimate an absolute confidence $\delta t={\rm{max}}|\bar{t}-\bar t_{95\%}|$.

\subsection{Statistical models for the PSO time data}
\subsubsection{PSO time data in the range $R_3$}

We begin by analyzing our PSO method. We looked for solutions to $x^3+y^3+z^3=2$ in the lowest range $R_3$. In this case, the method was tested $N=10^4$ times. The produced set of running times $\{t_i\}_{i=1}^{N}$ is divided into bins with width $\Delta t = 0.00019\,s$ and its PDF histogram is shown in figure \ref{fig:DFO-Hist1}. 

\begin{figure}[H]
 	\centering
 	\begin{subfigure}[b]{0.47\textwidth}
 		\includegraphics[width=\textwidth]{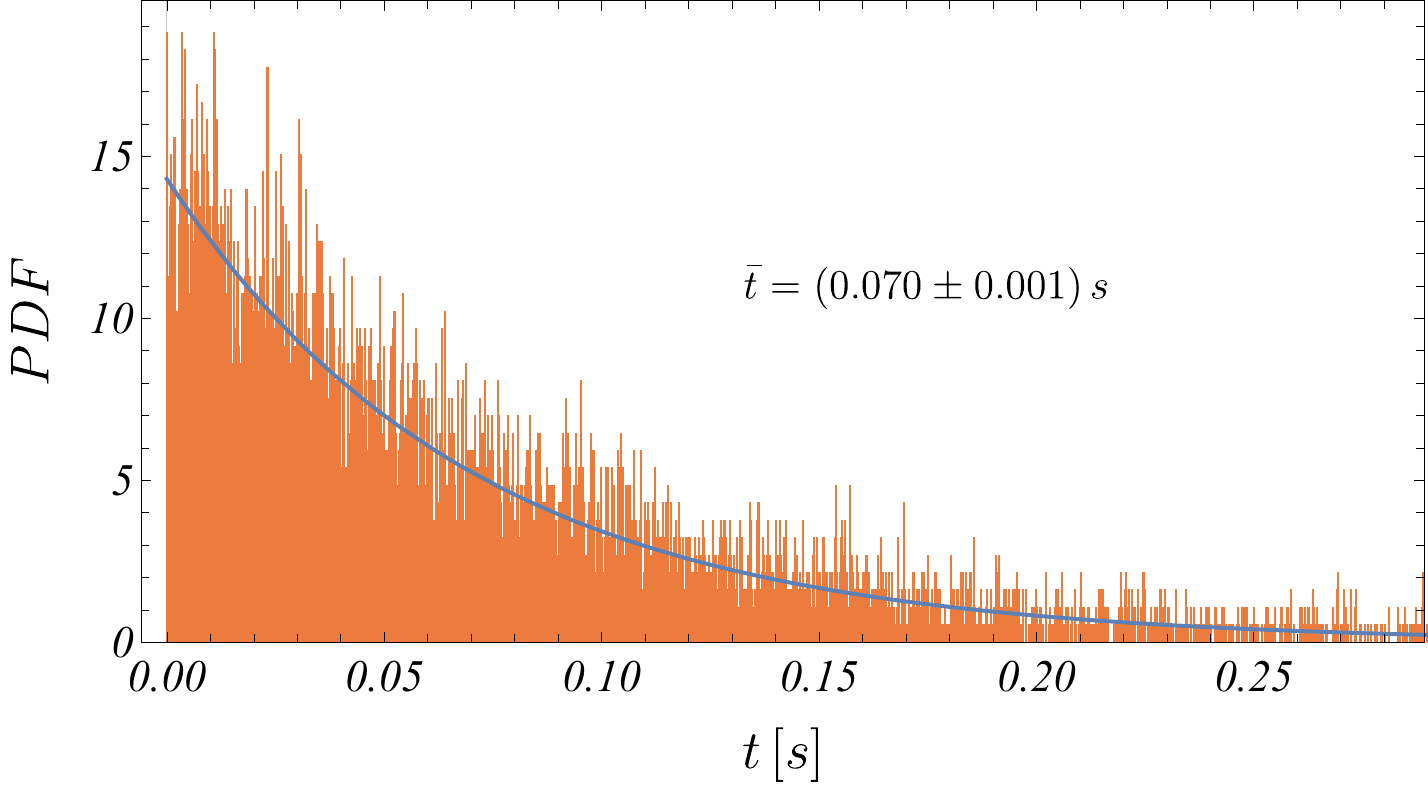}
 		\caption{PSO time data exponential fit.}
 		\label{fig:DFO-Exp1}
 	\end{subfigure}
 	\hspace{0.8 mm}
 	\begin{subfigure}[b]{0.47\textwidth}
 		\includegraphics[width=\textwidth]{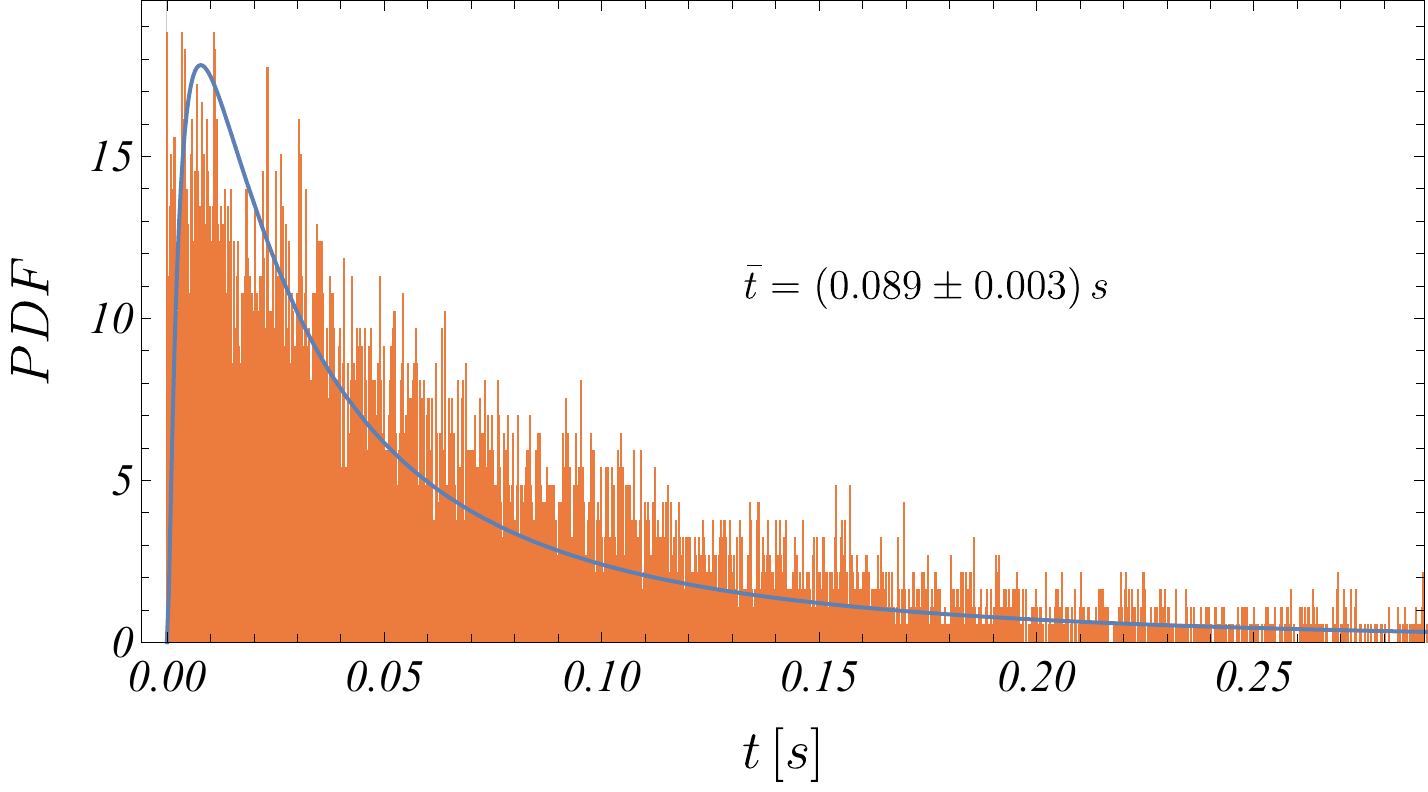}
 		\caption{PSO time data log-normal fit.}
 		\label{fig:DFO-LogN1}
 	\end{subfigure}
 	\caption{Statistical models for the PSO time data in the range $R_3$. (a) shows a simple one-parameter exponential probability distribution $f(t;\lambda)$ with $\lambda=14.301\, s^{-1}$ and mean time for finding a solution $\bar{t}=0.070 \,s$. (b) shows a two-parameter log-normal distribution fit $f(t;\alpha,\beta)$ with $\alpha=-3.229,\,\beta=1.275$, and mean time for finding a solution $\bar{t}=0.089 \,s$.}\label{fig:DFO1}
 \end{figure}

The two statistical models, describing the PSO time data, are shown in figure \ref{fig:DFO1}. The first one is a simple one-parameter exponential model $f(t;\lambda)$ with $\lambda=14.301 \,s^{-1}$ (figure \ref{fig:DFO-Exp1}). The second one is a two-parameter log-normal distribution $f(t;\alpha,\beta)$ with $\alpha=-3.229,\,\beta=1.275$ (figure \ref{fig:DFO-LogN1}).

{\bf{Exponential model}}

The expected value, the variance and the median of an exponentially distributed random variable $t$ with rate parameter $\lambda=14.301 \, s^{-1}$ yield
\begin{equation}
\bar t={\mathop{\rm E}\nolimits} [t] = \frac{1}{\lambda } = 0.070\,s,\quad {\mathop{\rm Var}\nolimits} [t] = \frac{1}{{{\lambda ^2}}} = 0.005\,{s^2},\quad {\mathop{\rm Med}\nolimits}[t] = \frac{{\ln 2}}{\lambda } = 0.048\,s.
\end{equation}

The $95\% $ confidence intervals for $\lambda$ and $\bar t$ are given by  
\begin{equation}
{\lambda _{lower}} \le \lambda  \le {\lambda _{upper}},\quad {{{\lambda^{-1} _{upper}}}} \le \bar t \le {{{\lambda^{-1} _{lower}}}},
\end{equation}
where
\begin{equation}
{\lambda _{lower}} = \lambda \left( {1 - \frac{{1.96}}{{\sqrt {N} }}} \right) \approx 14.021\,{s^{ - 1}},\quad {\lambda _{upper}} = \lambda \left( {1 + \frac{{1.96}}{{\sqrt {N} }}} \right) \approx 14.581\,{s^{ - 1}},
\end{equation}
\begin{equation}
{{{\lambda^{-1} _{upper}}}} = 0.069\,s,\quad {{{\lambda^{-1} _{lower}}}} = 0.071\,s,
\end{equation}
with absolute confidences $\delta \lambda=|\lambda-\lambda_{upper}|=0.280 \,s^{-1}$ and $\delta t=|\bar t-\lambda_{lower}^{-1}|=0.001\,s$. Therefore, we can write our results for the Poisson distributed PSO time data as
\begin{equation}
\lambda  = (14.301 \pm 0.280)\,{s^{ - 1}},\quad \bar t = (0.070\, \pm 0.001)\,s.
\end{equation}

Because $\lambda$ and $\bar t$ are inversely proportional to each other, from now on we will be interested only in $\bar t$. 

{\bf{Log-normal model}}

The two-parameter log-normal distribution $f(t;\alpha,\beta)$ for the PSO time data has estimated parameters $\alpha=-3.229,\,\beta=1.275$, which is depicted in figure \ref{fig:DFO-LogN1}. The parameters can be obtained from the sample data via
\begin{equation}
\alpha  = \frac{1}{N}\sum\limits_{i = 1}^{N} {\ln {t_i}} =  - 3.229,\quad \beta  = \frac{1}{{\sqrt {N} }}\sqrt {\sum\limits_{i = 1}^{N} {{{(\ln {t_i} - \alpha )}^2}} } = 1.275.
\end{equation}

Consequently, the mean time $\bar{t}$ for finding a solution and its standard deviation are given by
\begin{equation}
\bar t = {\mathop{\rm E}\nolimits} [t] = {e^{\alpha  + \frac{{{\beta ^2}}}{2}}} = 0.089\,s,\quad {\mathop{\rm SD}\nolimits} [t]=\sqrt{{\mathop{\rm Var}\nolimits} [t]} = {e^{\alpha  + \frac{\beta ^2}{2}}}\sqrt{{e^{{\beta ^2}}} - 1} = 0.180\,{s}.
\end{equation}

Furthermore, the median and the mode yield
\begin{equation}
{\mathop{\rm Med}\nolimits} [t] = {e^\alpha } = 0.040\,s,\quad {\mathop{\rm Mode}\nolimits} [t] = {e^{\alpha  - {\beta ^2}}} = 0.008\,s,
\end{equation}

The standard scatter intervals for the log-normal distribution are given by 	
\begin{eqnarray}
{t_{68\% }} \in [{e^{\alpha  - \beta }},{e^{\alpha  + \beta }}] = [0.011\,s,0.142\,s],
\end{eqnarray}
for the $68\%$ confidence interval, and
\begin{eqnarray}
{t_{95\% }} \in [{e^{\alpha  - 2\beta }},{e^{\alpha  + 2\beta }}] = [0.003\,s,0.507\,s],
\end{eqnarray}
for the $95 \%$ confidence interval. 

The efficient $95\%$ confidence interval for $\bar t$ yields
\begin{equation}\label{eqCoxConfInt1}
{\bar t_{95\% }} \in {e^{\left[ {\alpha  + \frac{{{\beta ^2}}}{2} - 1.96\sqrt {\frac{{{\beta ^2}}}{{{N}}} + \frac{{{\beta ^4}}}{{2({N} - 1)}}} ,{\kern 1pt} \alpha  + \frac{{{\beta ^2}}}{2} + 1.96\sqrt {\frac{{{\beta ^2}}}{{{N}}} + \frac{{{\beta ^4}}}{{2({N} - 1)}}} } \right]}} = [0.086\,s,\,0.092\,s],
\end{equation}
with absolute confidence $\delta t={\rm{max}}|\bar{t}-\bar t_{95\%}|=0.003\, s$, thus
\begin{equation}
\bar t = (0.089\, \pm 0.003)\,s.
\end{equation}	
The relevant data is collected in table \ref{table:1}.
\begin{table}[h]
	\small\centering
	\begin{tabular}{ ||c| c| c| c| c| c| c|c||} 
		\hline
		\backslashbox{Dist.}{Param.} &$\alpha$  &$\beta$  &$\lambda\ [s^{-1}]$ &$\bar t\pm\delta t\ [s]$ &${\bar t_{95\% }}\ [s]$& ${\mathop{\rm Med}\nolimits}[t] \ [s]$ &${\mathop{\rm Mode}\nolimits}[t]\ [s]$\\	\hline
		Exponential & -- & -- & 14.301& $0.070\, \pm 0.001$& [0.069,0.071]& 0.048& --\\	\hline
		Log-normal &-3.229 &1.275 &--& $0.089\, \pm 0.003$& [0.086,0.092]&0.040&0.008\\	\hline
	\end{tabular}
\caption{The relevant characteristics of the models for PSO in $R_3$.}
\label{table:1}
\end{table}
\subsubsection{PSO time data in the range $R_4$}

\begin{figure}[H]
	\centering
	\begin{subfigure}[b]{0.47\textwidth}
		\includegraphics[width=\textwidth]{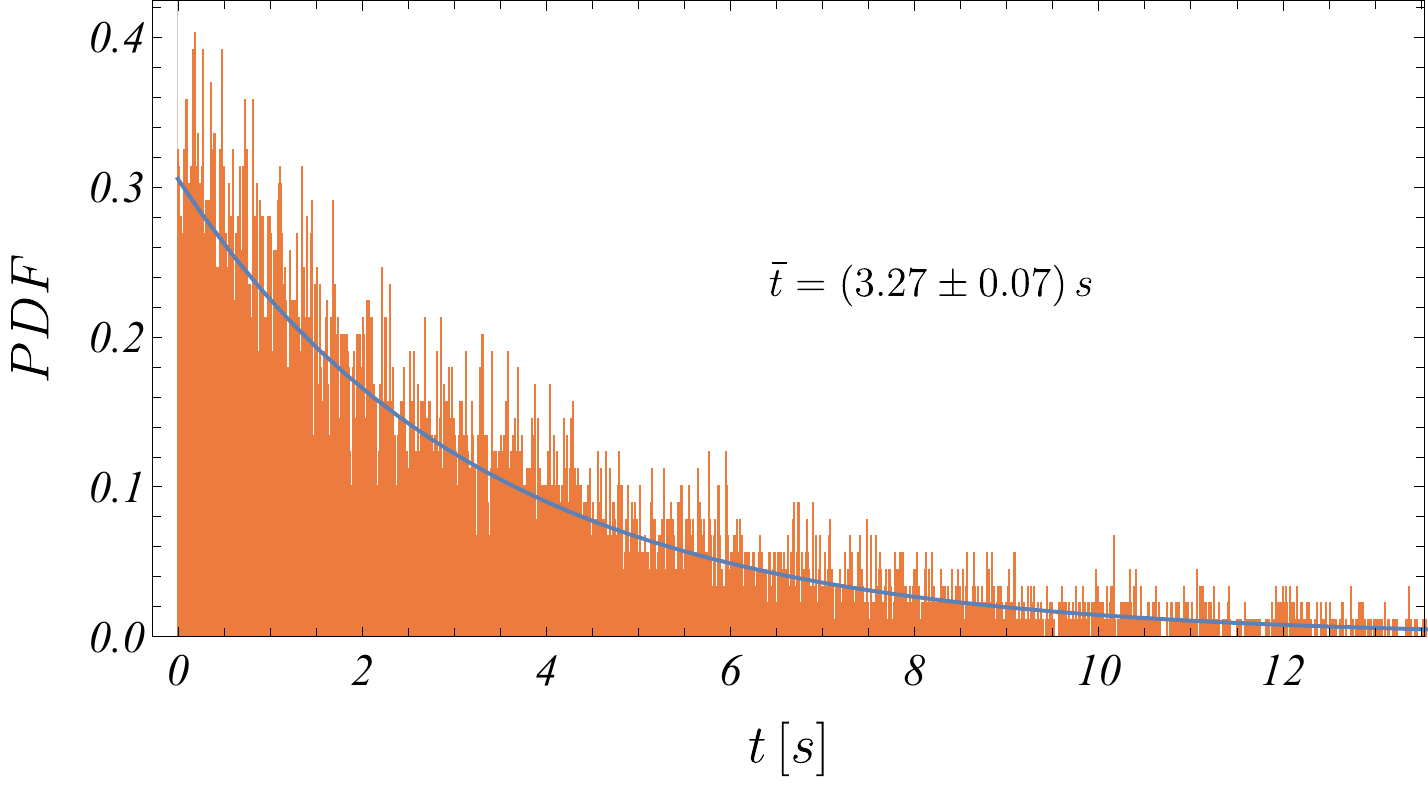}
		\caption{PSO time data exponential fit.}
		\label{fig:DFO-Exp2}
	\end{subfigure}
	\hspace{0.8 mm}
	\begin{subfigure}[b]{0.47\textwidth}
		\includegraphics[width=\textwidth]{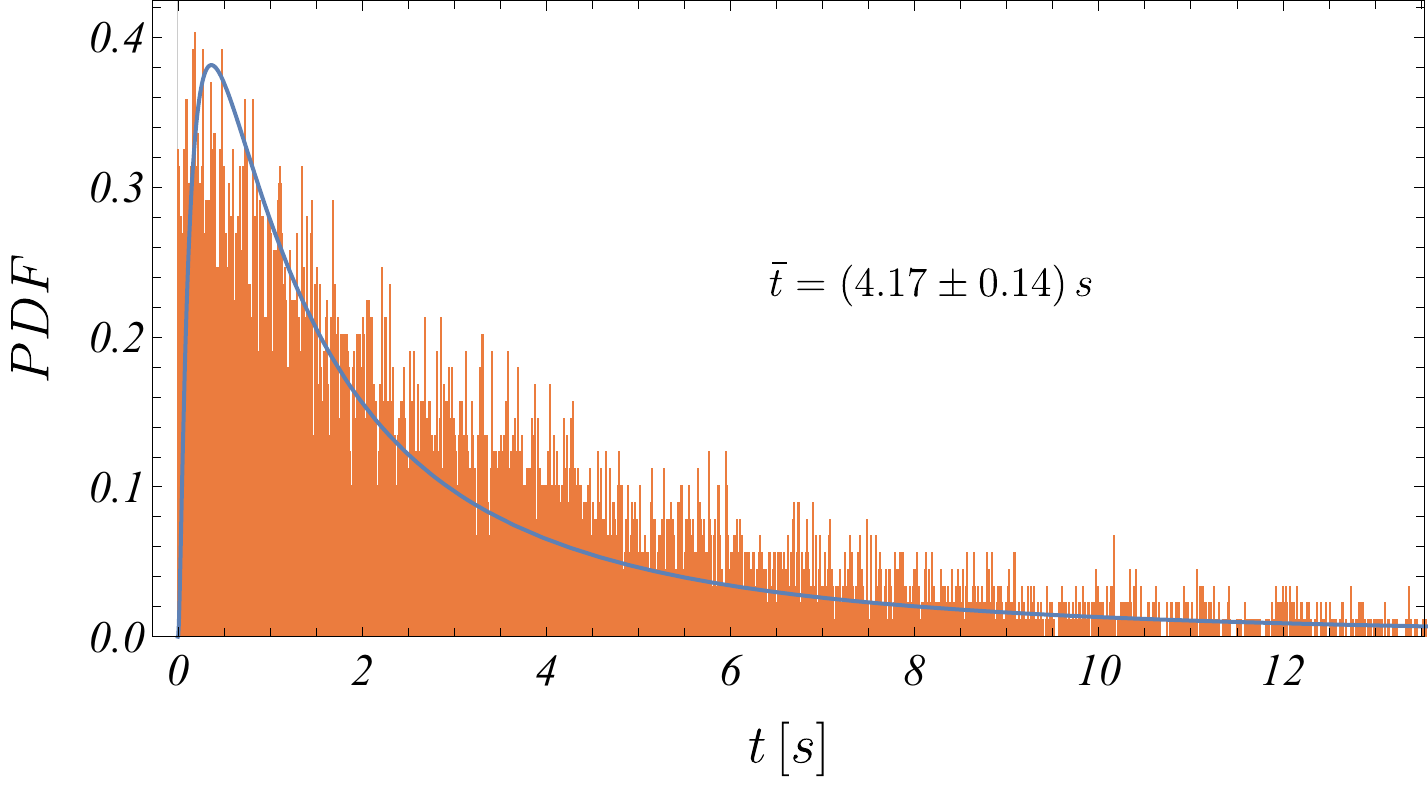}
		\caption{PSO time data log-normal fit.}
		\label{fig:DFO-LogN2}
	\end{subfigure}
	\caption{Statistical models for the PSO time data in the range $R_4$. (a) shows a one-parameter exponential model fit $f(t;\lambda)$ with $\lambda=0.31\, s^{-1}$ and mean time for finding a solution $\bar{t}=3.27 \,s$. (b) shows a two-parameter log-normal distribution fit $f(t;\alpha,\beta)$ with $\alpha=0.62,\,\beta=1.28$, and mean time for finding a solution $\bar{t}=4.17 \,s$.}\label{fig:DFO2}
\end{figure}

Next, we analyze the PSO time data, accumulated when looking for integer solutions to $x^3+y^3+z^3=2$ in the mid range $R_4$. In this case, the method was tested $N=10^4$ times. The produced time set $\{t_i\}_{i=1}^{N}$ is divided into bins with width $\Delta t = 0.01\,s$ and its PDF histogram is shown in figure \ref{fig:DFO-Hist2}. As in the previous case, we model the PSO time data by an exponential model, shown in figure \ref{fig:DFO-Exp2}, and a log-normal model, shown in figure \ref{fig:DFO-LogN2}, with the relevant characteristics collected in table \ref{table:2}.

\begin{table}[h]
	\small\centering
	\begin{tabular}{ ||c| c| c| c| c| c| c|c||} 
		\hline
		\backslashbox{Dist.}{Param.} &$\alpha$  &$\beta$  &$\lambda\ [s^{-1}]$ &$\bar t\pm\delta t \ [s]$ &${\bar t_{95\% }}\ [s]$& ${\mathop{\rm Med}\nolimits}[t] \ [s]$ &${\mathop{\rm Mode}\nolimits}[t] \ [s]$\\	\hline
		Exponential & -- & -- & 0.31& $3.27\, \pm 0.07$& [3.21,3.34]& 2.27& --\\	\hline
		Log-normal &0.62 &1.28 &--& $4.17\, \pm 0.14$& [4.03,4.31]&1.85&0.36\\	\hline
	\end{tabular}
	\caption{The relevant characteristics of the models for PSO in $R_4$.}
	\label{table:2}
\end{table}

\subsubsection{PSO time data in the range $R_5$}

\begin{figure}[H]
	\centering
	\begin{subfigure}[b]{0.47\textwidth}
		\includegraphics[width=\textwidth]{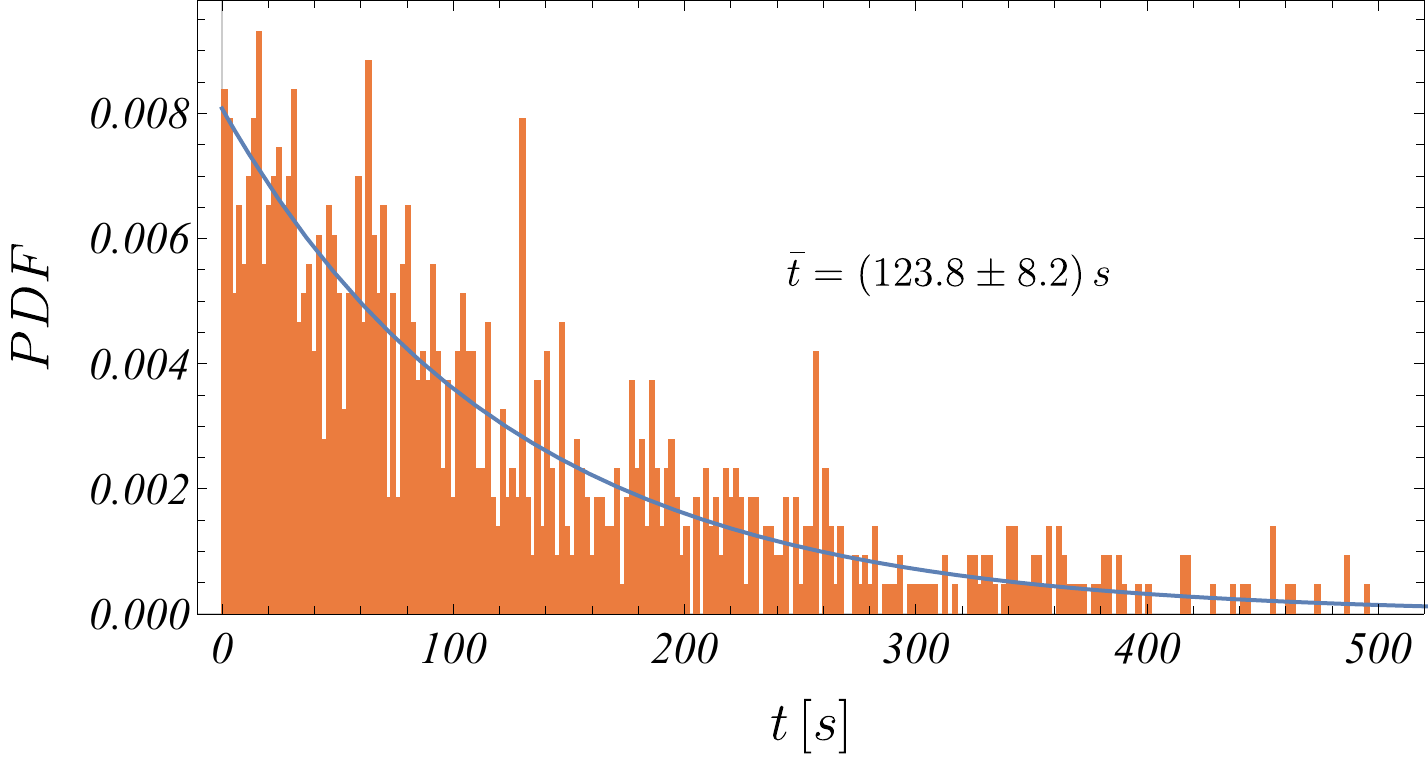}
		\caption{PSO time data exponential fit.}
		\label{fig:DFO-Exp3}
	\end{subfigure}
	\hspace{0.8 mm}
	\begin{subfigure}[b]{0.47\textwidth}
		\includegraphics[width=\textwidth]{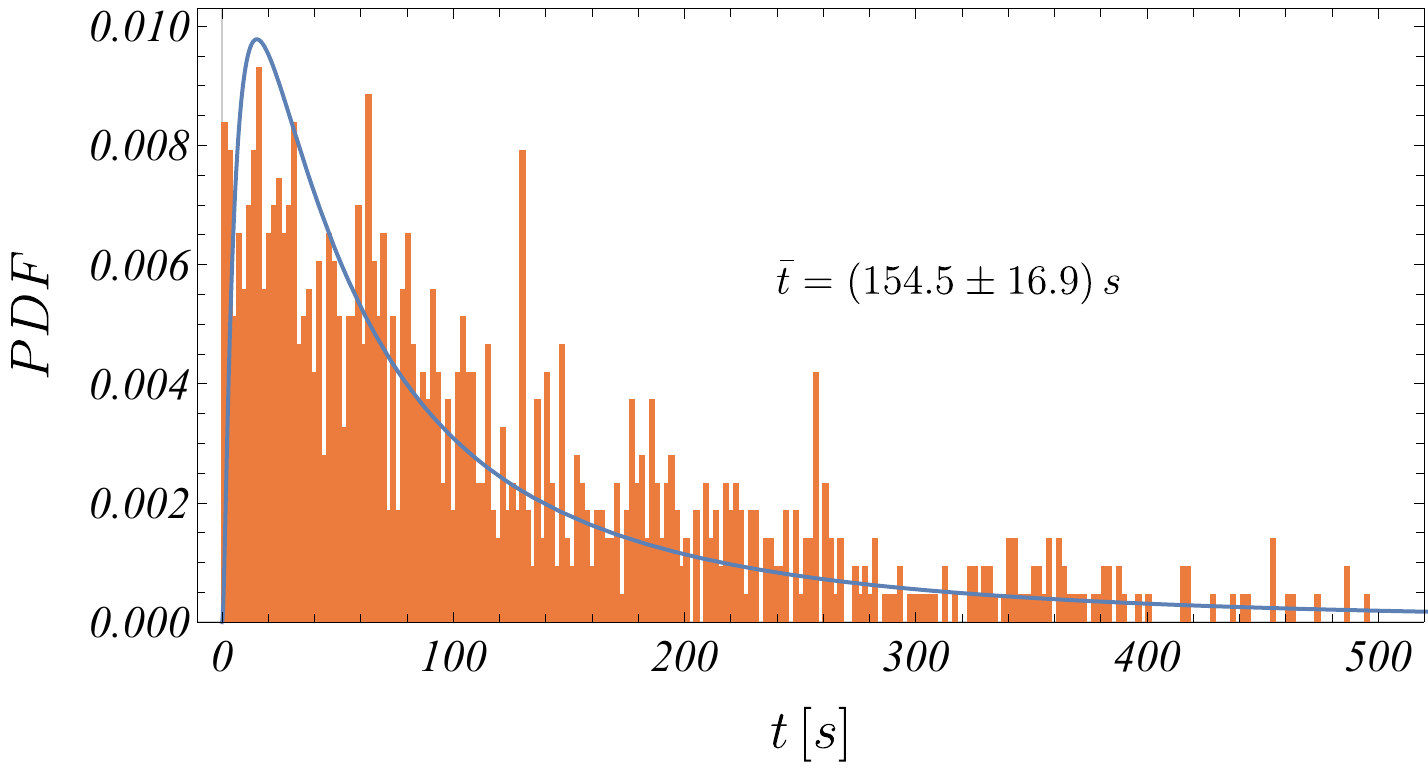}
		\caption{PSO time data log-normal fit.}
		\label{fig:DFO-LogN3}
	\end{subfigure}
	\caption{Statistical models for the PSO time data the range $R_5$. (a) shows a one-parameter exponential model fit $f(t;\lambda)$ with $\lambda=0.008\, s^{-1}$ and mean time for finding a solution $\bar{t}=123.8 \,s$. (b) shows a two-parameter log-normal distribution fit $f(t;\alpha,\beta)$ with $\alpha=4.27,\,\beta=1.25$, and mean time for finding a solution $\bar{t}=154.5 \,s$.}\label{fig:DFO3}
\end{figure}

We continue our analysis by looking for solutions to $x^3+y^3+z^3=2$ in the range $R_5$. In this case, the method was tested $N=10^3$ times. The produced time set $\{t_i\}_{i=1}^{N}$ is divided into bins with width $\Delta t = 2.15\,s$ and its PDF histogram is shown in figure \ref{fig:DFO-Hist3}. The statistical models, describing the PSO time data, are shown in figures \ref{fig:DFO-Exp3} and \ref{fig:DFO-LogN3} and their characteristics -- in table \ref{table:3}.

\begin{table}[h]
	\small\centering
	\begin{tabular}{ ||c| c| c| c| c| c| c|c||} 
		\hline
		\backslashbox{Dist.}{Param.} &$\alpha$  &$\beta$  &$\lambda\ [s^{-1}]$ &$\bar t\pm\delta t\ [s]$ &${\bar t_{95\% }}\ [s]$& ${\mathop{\rm Med}\nolimits}[t] \ [s]$ &${\mathop{\rm Mode}\nolimits}[t]\ [s]$\\	\hline
		Exponential & -- & -- & 0.008& $123.8\, \pm 8.2$& [116.6,132.0]& 85.8& --\\	\hline
		Log-normal &4.27 &1.25 &--& $154.5\, \pm 17.0$& [139.5,171.5]&71.1&15.1\\	\hline
	\end{tabular}
	\caption{The relevant characteristics of the models for PSO in $R_5$.}
	\label{table:3}
\end{table}
\subsection{Statistical models for the SA algorithm time data}

\subsubsection{SA time data in the range $R_3$}

In this section we focus on the time data accumulated from the SA algorithm (without restarts). The analysis mimics the one for the PSO method. 

\begin{figure}[H]
	\centering
	\begin{subfigure}[b]{0.47\textwidth}
		\includegraphics[width=\textwidth]{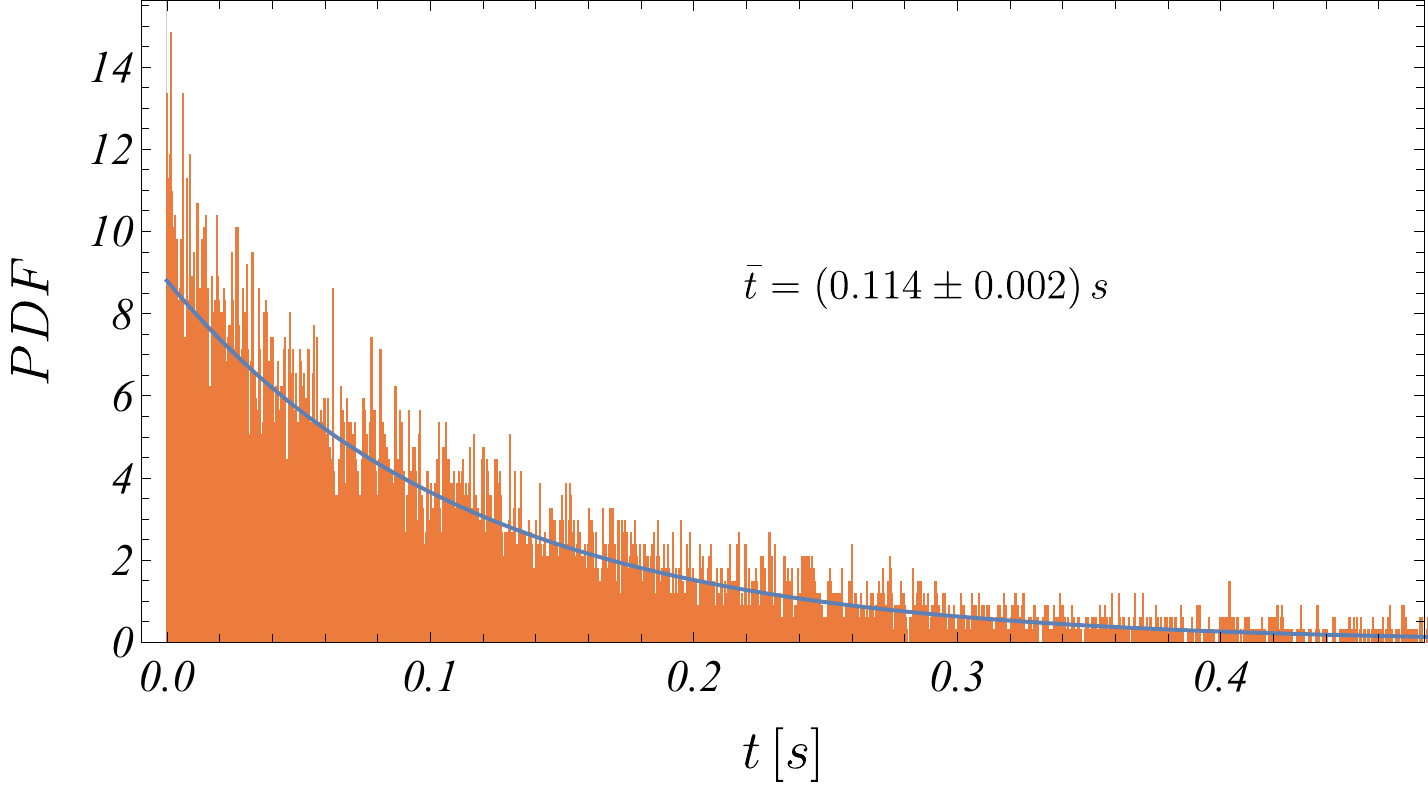}
		\caption{Exponential fit for SA.}
		\label{fig:SANORES-Exp1}
	\end{subfigure}
	\hspace{0.8 mm}
	\begin{subfigure}[b]{0.47\textwidth}
		\includegraphics[width=\textwidth]{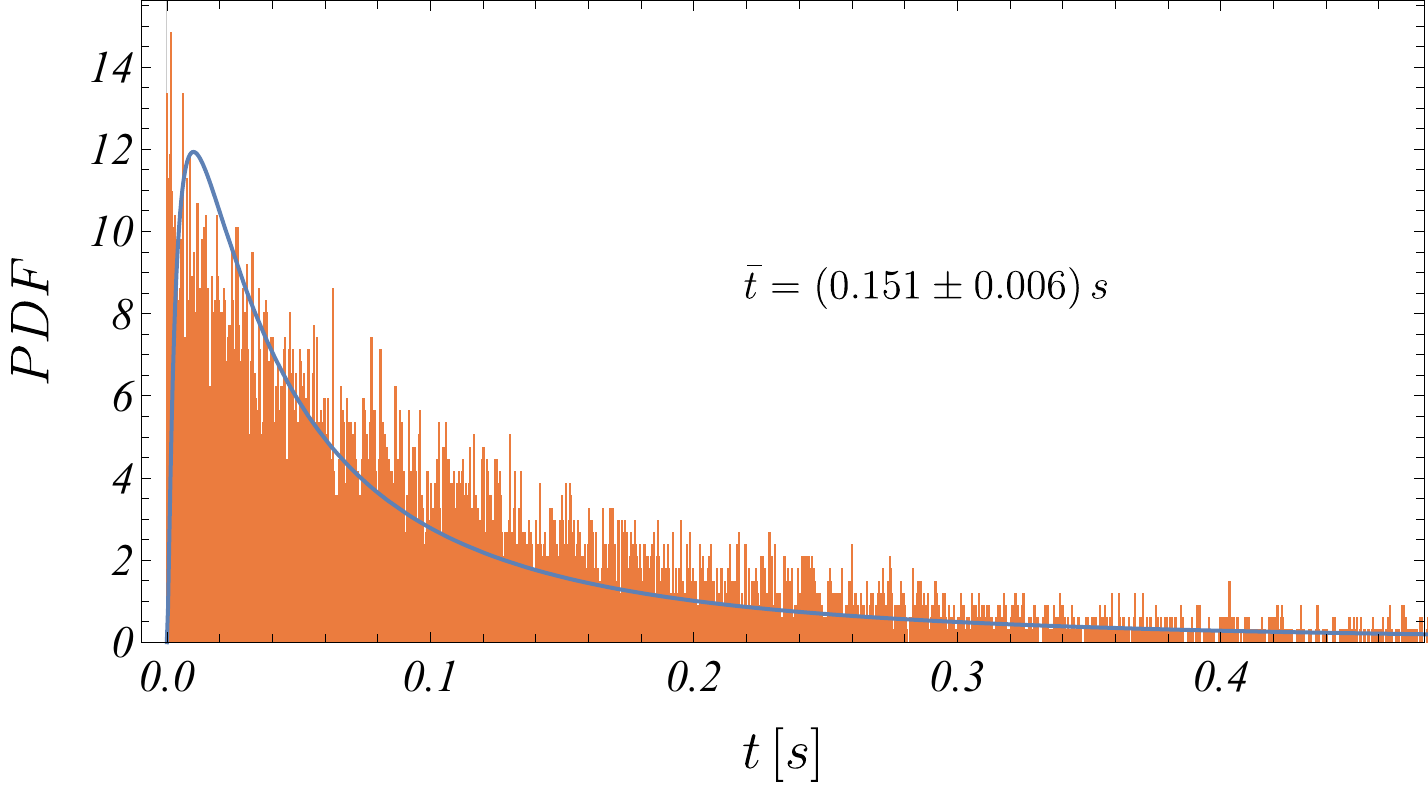}
		\caption{Log-normal fit for SA.}
		\label{fig:SANORES-LogN1}
	\end{subfigure}
	\caption{Statistical models for the SA time data in the range $R_{3}$. (a) shows an exponential model fit $f(t;\lambda)$ with $\lambda=8.800\, s^{-1}$ and mean time for finding a solution $\bar{t}=0.114 \,s$. (b) shows a  log-normal distribution fit $f(t;\alpha,\beta)$ with $\alpha=-2.791,\,\beta=1.343$, and mean time for finding a solution $\bar{t}=0.151 \,s$.}\label{fig:SANORES1}
\end{figure}

In the lowest range $R_{3}$ the produced time data $\{t_i\}_{i=1}^{N}$ is divided into bins with width $\Delta t = 0.00034\,s$ and its PDF histogram is shown in figure \ref{fig:SANORES-Hist1}. The chosen statistical models, describing the time data of the algorithm, are shown in figures \ref{fig:SANORES-Exp1} and \ref{fig:SANORES-LogN1}, correspondingly. Their characteristics are collected in table \ref{table:4}.

\begin{table}[h]
	\small\centering
	\begin{tabular}{ ||c| c| c| c| c| c| c|c||} 
		\hline
		\backslashbox{Dist.}{Param.} &$\alpha$  &$\beta$  &$\lambda\ [s^{-1}]$ &$\bar t\pm\delta t\ [s]$ &${\bar t_{95\% }}\ [s]$& ${\mathop{\rm Med}\nolimits}[t]\  [s]$ &${\mathop{\rm Mode}\nolimits}[t]\ [s]$\\	\hline
		Exponential & -- & -- & 8.800& $0.114\, \pm 0.002$& [0.111,0.116]& 0.079& --\\	\hline
		Log-normal &-2.791 &1.343 &--& $0.151\, \pm 0.006$& [0.146,0.157]&0.061&0.010\\	\hline
	\end{tabular}
	\caption{The relevant characteristics of the models for SA in $R_3$.}
	\label{table:4}
\end{table}

\subsubsection{SA time data in the range $R_4$}

\begin{figure}[H]
	\centering
	\begin{subfigure}[b]{0.47\textwidth}
		\includegraphics[width=\textwidth]{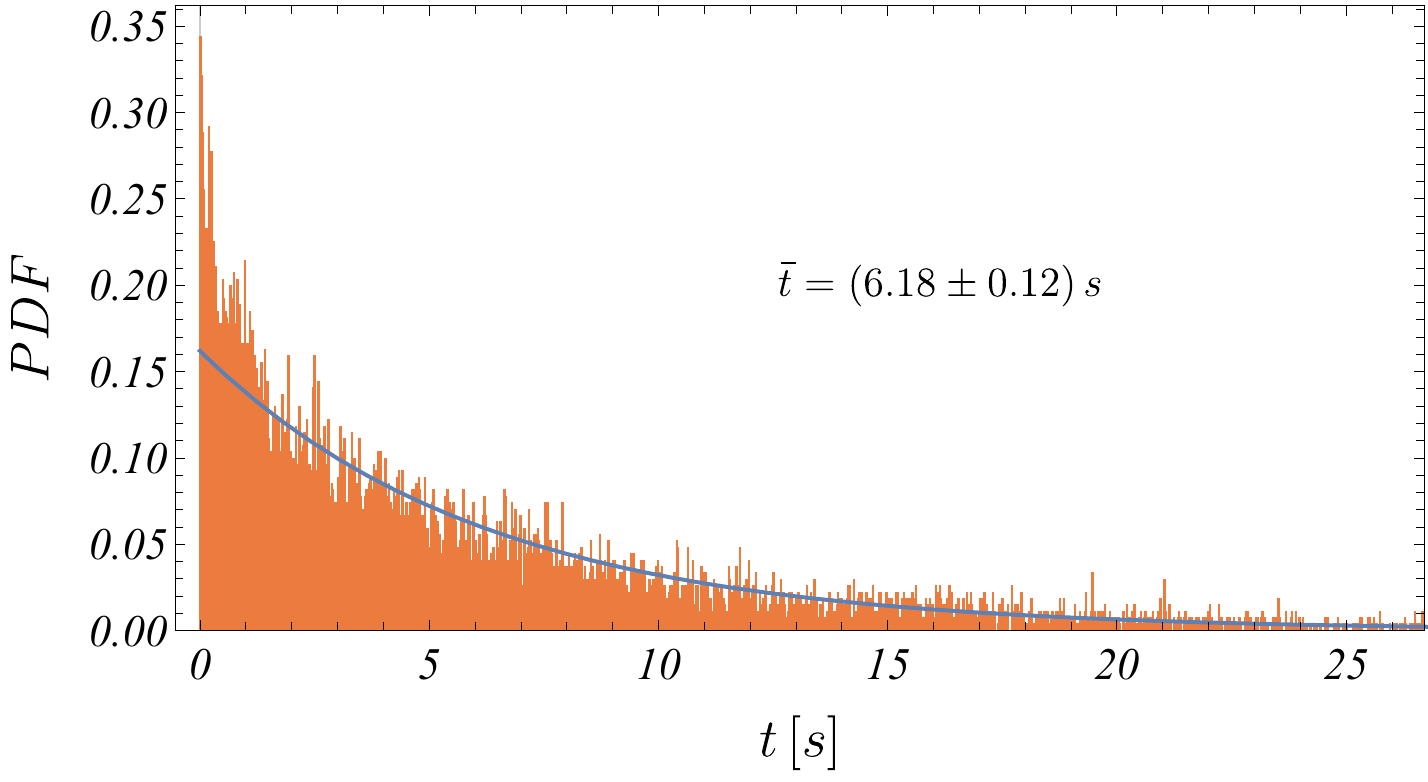}
		\caption{Exponential fit for SA.}
		\label{fig:SANORES-Exp2}
	\end{subfigure}
	\hspace{0.8 mm}
	\begin{subfigure}[b]{0.47\textwidth}
		\includegraphics[width=\textwidth]{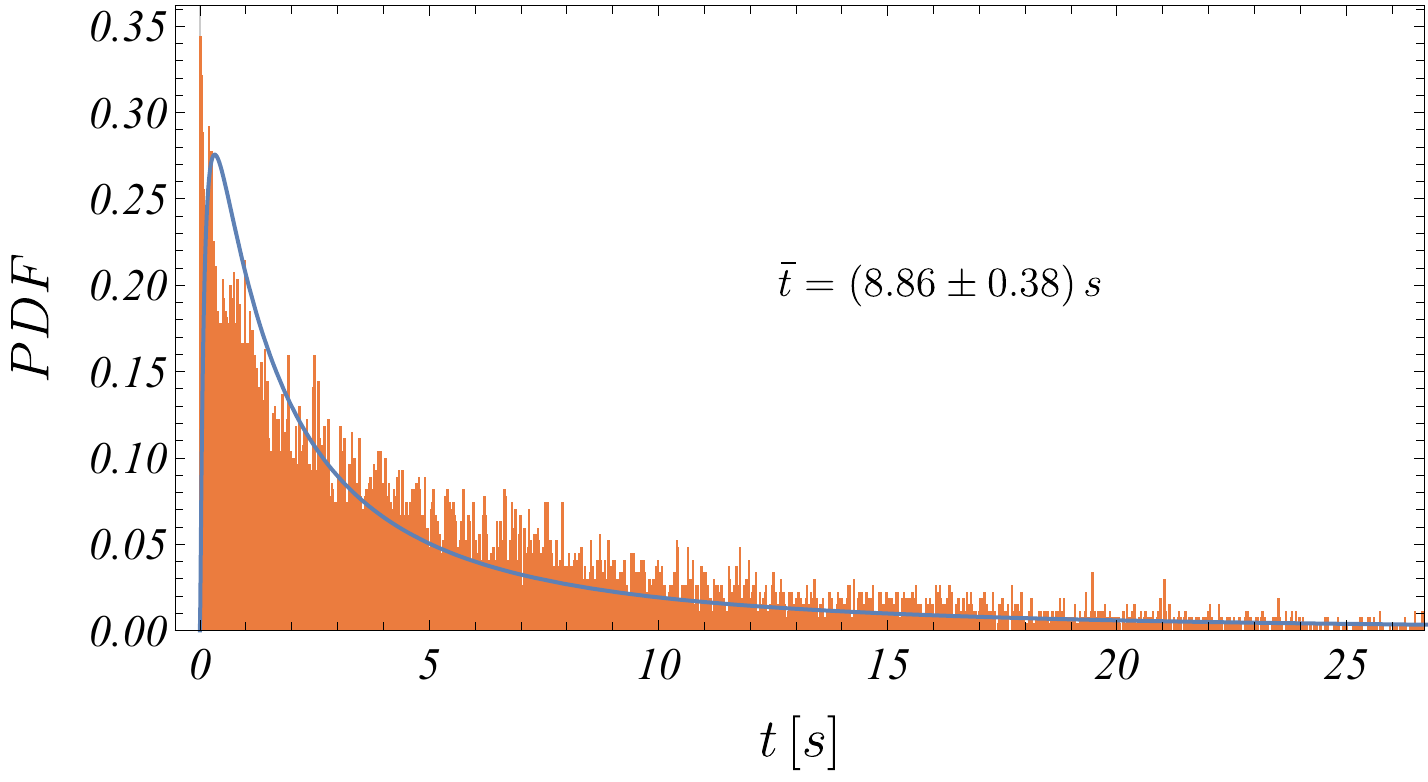}
		\caption{Log-normal fit for SA.}
		\label{fig:SANORES-LogN2}
	\end{subfigure}
	\caption{Statistical models for the SA time data in the range $R_4$. (a) shows an exponential model fit $f(t;\lambda)$ with $\lambda=0.16\, s^{-1}$ and mean time for finding a solution $\bar{t}=6.18\,s$. (b) shows a log-normal distribution fit $f(t;\alpha,\beta)$ with $\alpha=1.08,\,\beta=1.49$, and mean time for finding a solution $\bar{t}=8.86 \,s$.}\label{fig:SANORES2}
\end{figure}

We consider the range $R_4$ with $N=10^4$. The produced time set $\{t_i\}_{i=1}^{N}$ is divided into bins with width $\Delta t = 0.027\,s$ and its PDF histogram is shown in figure \ref{fig:SANORES-Hist2}. The considered statistical models, describing the SA time data, are shown in figures \ref{fig:SANORES-Exp2} and \ref{fig:SANORES-LogN2} and their characteristics -- in table \ref{table:5}.

\begin{table}[h]
	\small\centering
	\begin{tabular}{ ||c| c| c| c| c| c| c|c||} 
		\hline
		\backslashbox{Dist.}{Param.} &$\alpha$  &$\beta$  &$\lambda\ [s^{-1}]$ &$\bar t\pm\delta t\ [s]$ &${\bar t_{95\% }}\ [s]$& ${\mathop{\rm Med}\nolimits}[t] \ [s]$ &${\mathop{\rm Mode}\nolimits}[t]\ [s]$\\	\hline
		Exponential & -- & -- & 0.16& $6.18\, \pm 0.12$& [6.06,6.30]& 4.28& --\\	\hline
		Log-normal &1.08 &1.49 &--& $8.86\, \pm 0.38$& [8.50,9.25]&2.94&0.32\\	\hline
	\end{tabular}
	\caption{The relevant characteristics of the models for SA in $R_4$.}
	\label{table:5}
\end{table}

\subsubsection{SA time data in the range $R_5$}

\begin{figure}[H]
	\centering
	\begin{subfigure}[b]{0.47\textwidth}
		\includegraphics[width=\textwidth]{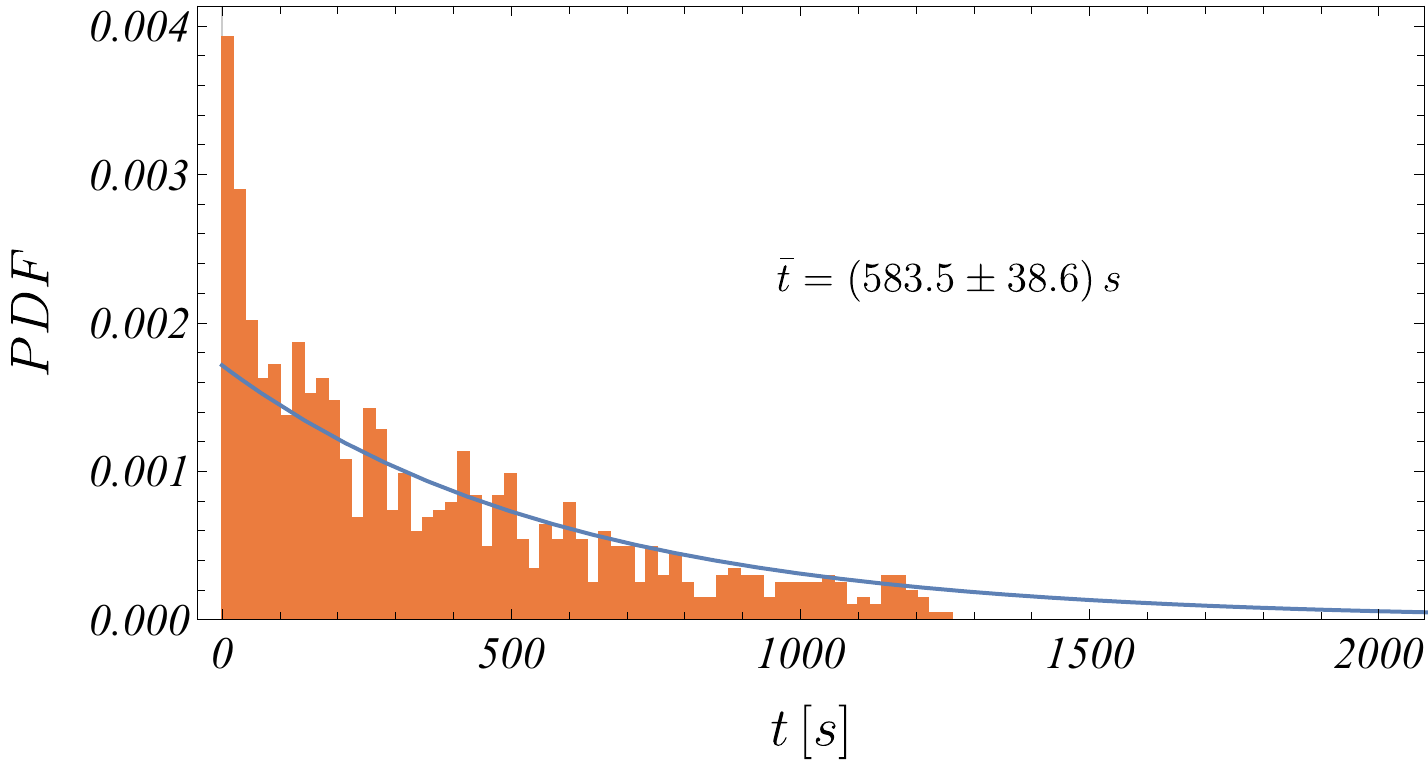}
		\caption{Exponential distribution fit for SA.}
		\label{fig:SANORES-Exp3}
	\end{subfigure}
	\hspace{0.8 mm}
	\begin{subfigure}[b]{0.47\textwidth}
		\includegraphics[width=\textwidth]{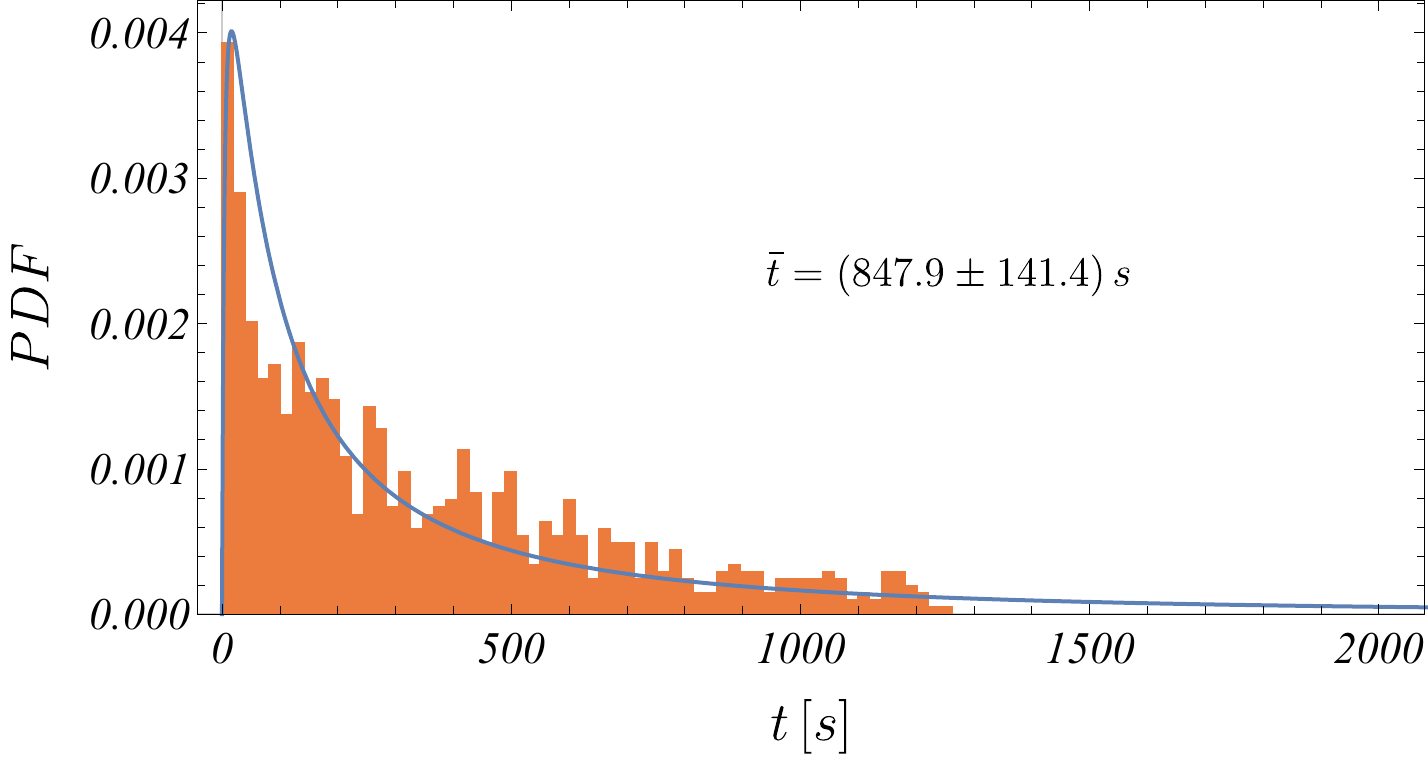}
		\caption{Log-normal distribution fit for SA.}
		\label{fig:SANORES-LogN3}
	\end{subfigure}
	\caption{Statistical models for the SA method in the range $R_5$. (a) shows an exponential model fit $f(t;\lambda)$ with $\lambda=0.0017\, s^{-1}$, and mean time for finding a solution $\bar{t}=583.5\,s$. (b) shows a log-normal distribution fit $f(t;\alpha,\beta)$ with $\alpha=5.43,\,\beta=1.62$, and mean time for finding a solution $\bar{t}=847.9 \,s$.}\label{fig:SANORES3}
\end{figure}

Next, we analyze the SA data in the range $R_5$ with $N=10^3$ tests. The produced time set $\{t_i\}_{i=1}^{N}$ is divided into bins with width $\Delta t = 20.4\,s$ and its PDF histogram is shown in figure \ref{fig:SANORES-Hist3}. The statistical models are shown in figures \ref{fig:SANORES-Exp3} and  \ref{fig:SANORES-LogN3}, correspondingly. Table \ref{table:6} shows their characteristics.

\begin{table}[h]
	\small\centering
	\begin{tabular}{ ||c| c| c| c| c| c| c|c||} 
		\hline
		\backslashbox{Dist.}{Param.} &$\alpha$  &$\beta$  &$\lambda\ [s^{-1}]$ &$\bar t\pm\delta t \ [s]$ &${\bar t_{95\% }}\ [s]$& ${\mathop{\rm Med}\nolimits}[t] \ [s]$ &${\mathop{\rm Mode}\nolimits}[t]\ [s]$\\	\hline
		Exponential & -- & -- & 0.0017& $583.5\, \pm 38.6$& [549.4,622.0]& 404.4& --\\	\hline
		Log-normal &5.43 &1.62 &--& $847.9\, \pm 141.4$& [728.6,989.3]&228.2&16.5\\	\hline
	\end{tabular}
	\caption{The relevant characteristics of the models for SA in $R_5$.}
	\label{table:6}
\end{table}

\subsection{Statistical models for the rSA algorithm time data}

\subsubsection{rSA time data in the range $R_3$}

\begin{figure}[H]
	\centering
	\begin{subfigure}[b]{0.47\textwidth}
		\includegraphics[width=\textwidth]{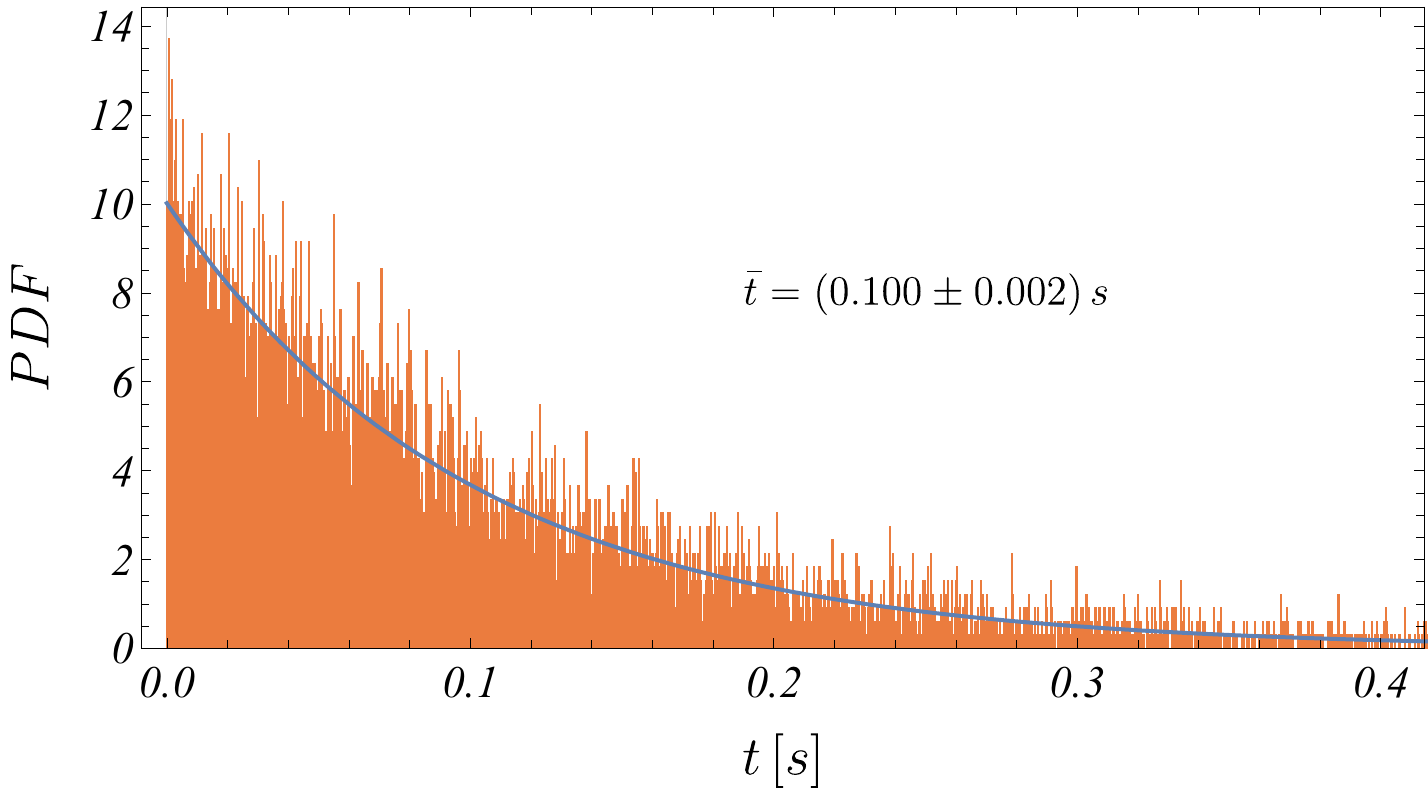}
		\caption{Exponential fit for rSA.}
		\label{fig:SARES-Exp1}
	\end{subfigure}
	\hspace{0.8 mm}
	\begin{subfigure}[b]{0.47\textwidth}
		\includegraphics[width=\textwidth]{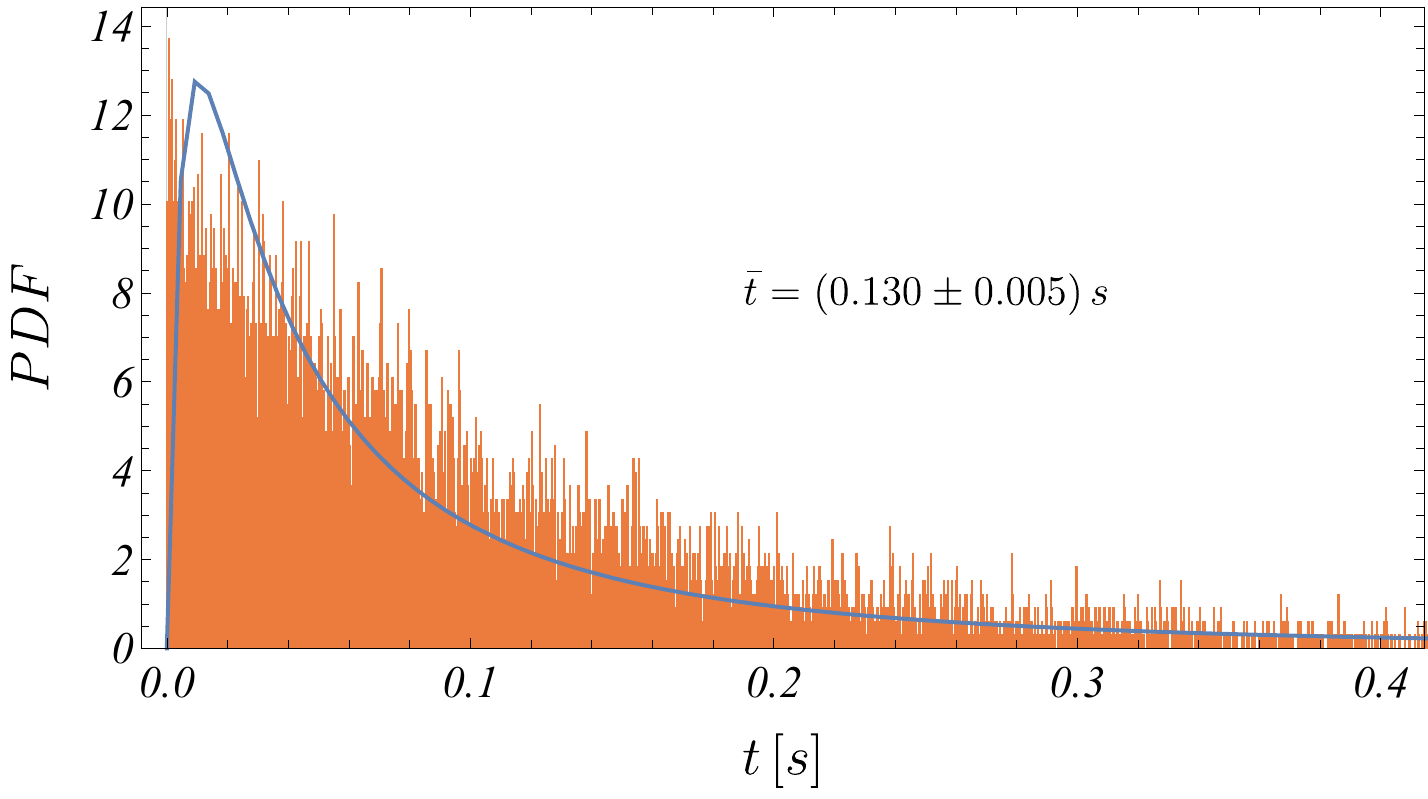}
		\caption{Log-normal fit for rSA.}
		\label{fig:SARES-LogN1}
	\end{subfigure}
	\caption{Statistical models for the rSA method in the range $R_3$. (a) shows an exponential model fit $f(t;\lambda)$ with $\lambda=10.014\, s^{-1}$ and mean time for finding a solution $\bar{t}=0.100 \,s$. (b) shows a log-normal distribution fit $f(t;\alpha,\beta)$ with $\alpha=-2.886,\,\beta=1.298$, and mean time for finding a solution $\bar{t}=0.130 \,s$.}\label{fig:SARES1}
\end{figure}

Here, we consider the rSA method in the lowest range $R_3$ and $N=10^4$. The produced time set $\{t_i\}_{i=1}^{N}$ is divided into bins with width $\Delta t = 0.00033\,s$ and its PDF histogram is shown in figure \ref{fig:SARES-Hist1}. The relevant statistical models are depicted in figures \ref{fig:SARES-Exp1} and \ref{fig:SARES-LogN1} with their characteristics shown in table \ref{table:7}.

\begin{table}[h]
	\small\centering
	\begin{tabular}{ ||c| c| c| c| c| c| c|c||} 
		\hline
		\backslashbox{Dist.}{Param.} &$\alpha$  &$\beta$  &$\lambda \ [s^{-1}]$ &$\bar t\pm\delta t\  [s]$ &${\bar t_{95\% }} \ [s]$& ${\mathop{\rm Med}\nolimits}[t] \ [s]$ &${\mathop{\rm Mode}\nolimits}[t]\ [s]$\\	\hline
		Exponential & -- & -- & 10.014& $0.100\, \pm 0.002$& [0.098,0.102]& 0.070& --\\	\hline
		Log-normal &-2.886 &1.298 &--& $0.130\, \pm 0.005$& [0.125,0.134]&0.056&0.010\\	\hline
	\end{tabular}
	\caption{The relevant characteristics of the models for rSA in $R_3$.}
	\label{table:7}
\end{table}

\subsubsection{rSA time data in the range $R_4$}

\begin{figure}[H]
	\centering
	\begin{subfigure}[b]{0.47\textwidth}
		\includegraphics[width=\textwidth]{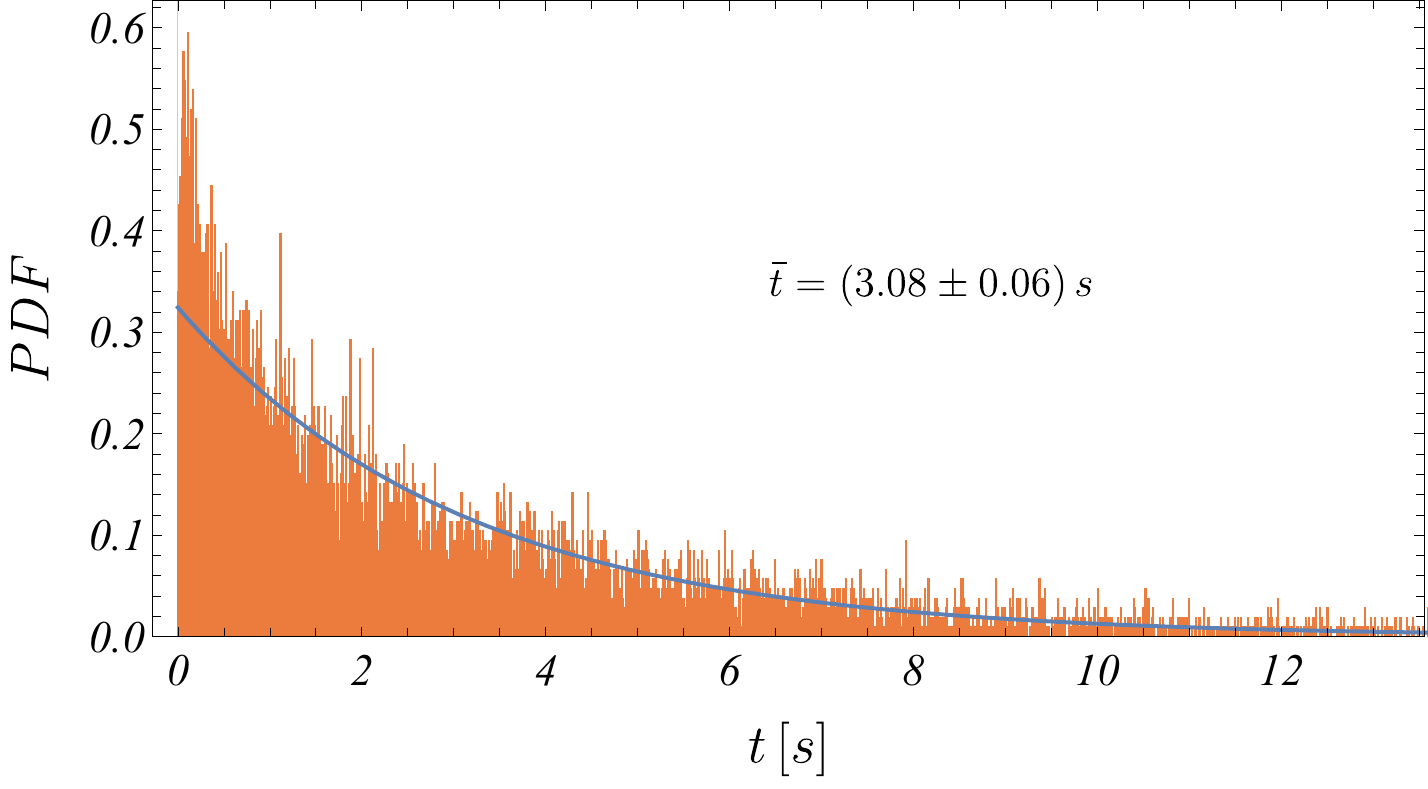}
		\caption{Exponential fit for rSA.}
		\label{fig:SARES-Exp2}
	\end{subfigure}
	\hspace{0.8 mm}
	\begin{subfigure}[b]{0.47\textwidth}
		\includegraphics[width=\textwidth]{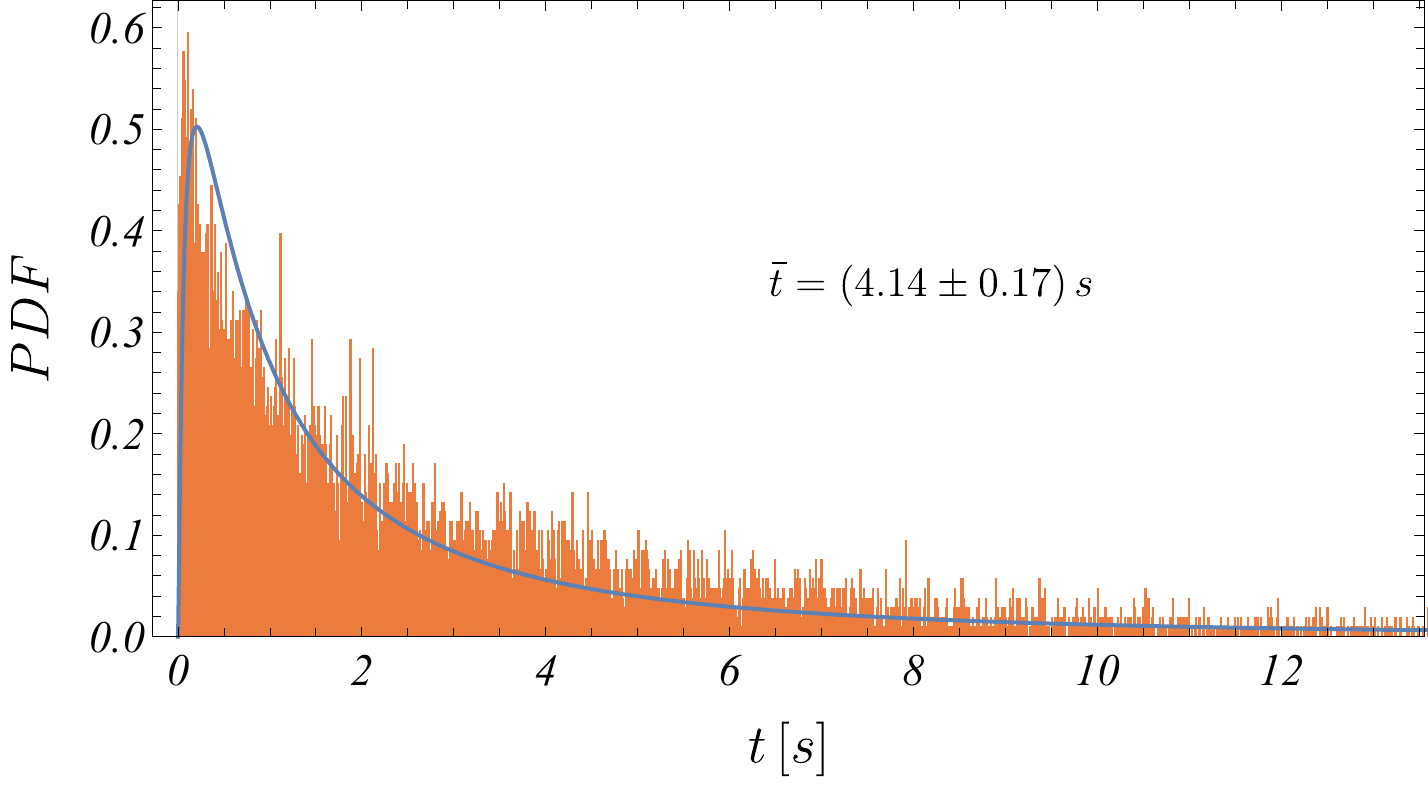}
		\caption{Log-normal fit for rSA.}
		\label{fig:SARES-LogN2}
	\end{subfigure}
	\caption{Statistical models for the rSA time data in the range $R_4$. (a) shows a one-parameter exponential model fit $f(t;\lambda)$ with $\lambda=0.32\, s^{-1}$ and mean time for finding a solution $\bar{t}=3.08\,s$. (b) shows a two-parameter log-normal distribution fit $f(t;\alpha,\beta)$ with $\alpha=0.42,\,\beta=1.41$, and mean time for finding a solution $\bar{t}=4.14 \,s$.}\label{fig:SARES2}
\end{figure}

We consider the range $R_4$ with $N=10^4$. The time data $\{t_i\}_{i=1}^{N}$ is divided into bins with width $\Delta t = 0.011\,s$ and its PDF histogram is shown in figure  \ref{fig:SARES-Hist2}. The statistical models are also shown in figures \ref{fig:SARES-Exp2} and \ref{fig:SARES-LogN2}. Table \ref{table:8} shows the relevant characteristics.

\begin{table}[h]
	\small\centering
	\begin{tabular}{ ||c| c| c| c| c| c| c|c||} 
		\hline
		\backslashbox{Dist.}{Param.} &$\alpha$  &$\beta$  &$\lambda\ [s^{-1}]$ &$\bar t\pm\delta t\ [s]$ &${\bar t_{95\% }}\ [s]$& ${\mathop{\rm Med}\nolimits}[t] \ [s]$ &${\mathop{\rm Mode}\nolimits}[t]\ [s]$\\	\hline
		Exponential & -- & -- & 0.32& $3.08\, \pm 0.06$& [3.02,3.14]& 2.14& --\\	\hline
		Log-normal &0.42 &1.41 &--& $4.14\, \pm 0.17$& [3.98,4.30]&1.52&0.21\\	\hline
	\end{tabular}
	\caption{The relevant characteristics of the models for rSA in $R_4$.}
	\label{table:8}
\end{table}

\subsubsection{rSA time data in the range $R_5$}

\begin{figure}[H]
	\centering
	\begin{subfigure}[b]{0.47\textwidth}
		\includegraphics[width=\textwidth]{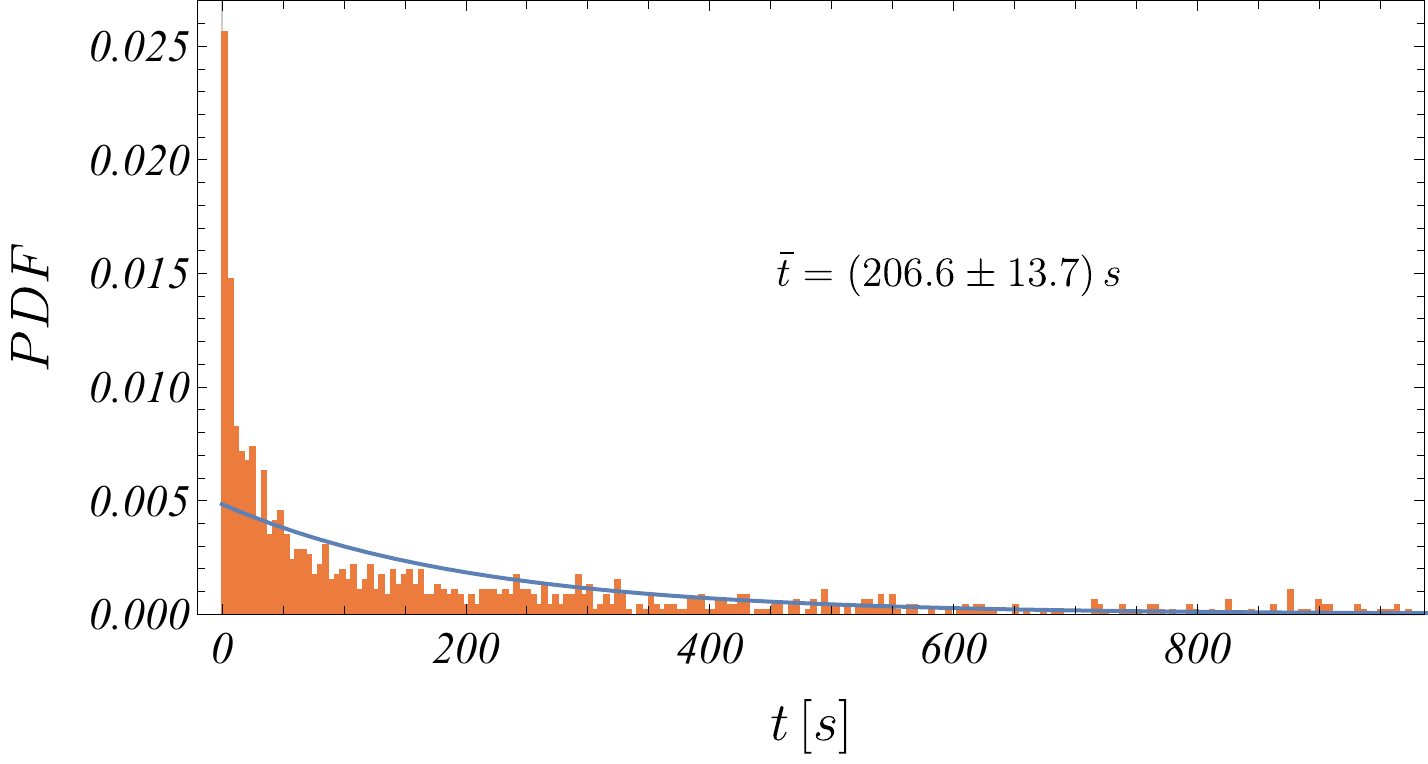}
		\caption{Exponential fit for rSA.}
		\label{fig:SARES-Exp3}
	\end{subfigure}
	\hspace{0.8 mm}
	\begin{subfigure}[b]{0.47\textwidth}
		\includegraphics[width=\textwidth]{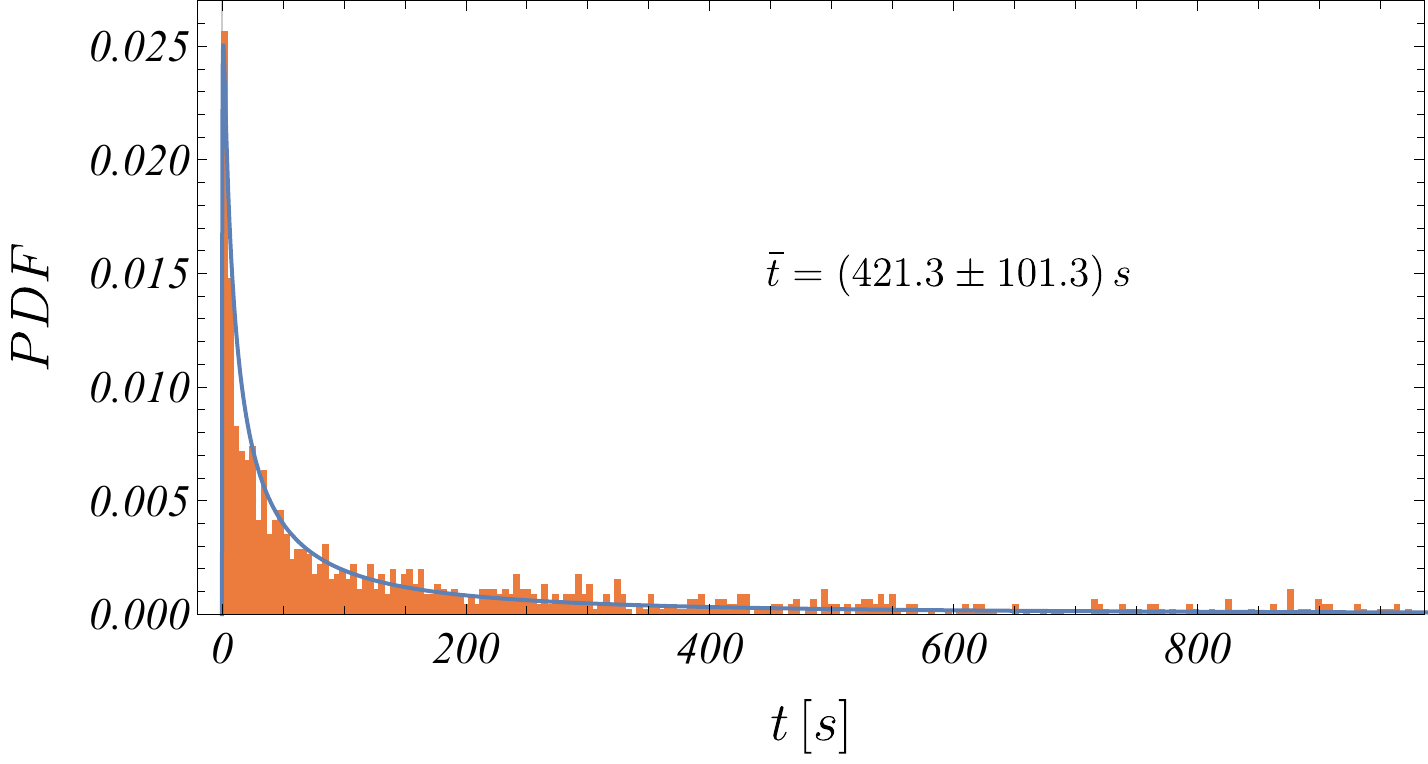}
		\caption{Log-normal fit for rSA.}
		\label{fig:SARES-LogN3}
	\end{subfigure}
	\caption{Statistical data models for rSA in the range $R_5$. (a) shows an exponential model fit $f(t;\lambda)$ with $\lambda=0.0048\, s^{-1}$ and mean time for finding a solution $\bar{t}=206.6\,s$. (b) shows a log-normal distribution fit $f(t;\alpha,\beta)$ with $\alpha=4.06,\,\beta=1.99$, and mean time for finding a solution $\bar{t}=421.3 \,s$.}\label{fig:SARES3}
\end{figure}

The final range is $R_5$ with $N=10^3$. The accumulated time data $\{t_i\}_{i=1}^{N}$ is divided into bins with width $\Delta t = 4.6\,s$ and its PDF histogram is shown in figure \ref{fig:SARES-Hist3}. 
The statistical data models are depicted in figures \ref{fig:SARES-Exp3} and \ref{fig:SARES-LogN3}. Their characteristics are in table \ref{table:9}.

\begin{table}[h!]
	\small\centering
	\begin{tabular}{ ||c| c| c| c| c| c| c|c||} 
		\hline
		\backslashbox{Dist.}{Param.} &$\alpha$  &$\beta$  &$\lambda\ [s^{-1}]$ &$\bar t\pm\delta t\  [s]$ &${\bar t_{95\% }}\ [s]$& ${\mathop{\rm Med}\nolimits}[t] \ [s]$ &${\mathop{\rm Mode}\nolimits}[t]\ [s]$\\	\hline
		Exponential & -- & -- & 0.0048& $206.6\, \pm 13.7$& [194.5,220.2]& 143.2& --\\	\hline
		Log-normal &4.06 &1.99 &--& $421.3\, \pm 101.3$& [341.0,522.6]&57.9&1.1\\	\hline
	\end{tabular}
	\caption{The relevant characteristics of the models for rSA in $R_5$.}
	\label{table:9}
\end{table}


 \newpage
\subsection{Fisher metric and model comparison}

In order to compare how dissimilar our statistical models are relative to each other, we need the explicit form of the Fisher information metric \cite{rao:1945,amari:2007,amari:2016,amari:2012, frieden:2010} for our distribution functions. Let $f(\vec{u};\vec\xi)$ be a PDF of some statistical model for a $d$-dimensional random variable $U$ with parameters $\vec \xi  = ({\xi ^1},{\xi ^2}, \ldots ,{\xi ^n})$. The Fisher metric is defined by the following integral over the range of $U$:
\begin{equation}\label{key}
{g_{ab}}(\vec \xi ) = \int\limits_U {\frac{{\partial \ln f(\vec u;\vec \xi )}}{{\partial {\xi ^a}}}} \frac{{\partial \ln f(\vec u;\vec \xi )}}{{\partial {\xi ^b}}}f(\vec u;\vec \xi ){\mathrm{d}^d}u\,,\quad a,b = 1, \ldots ,n.
\end{equation}

\noindent For one-dimensional models, consisting of a single free parameter, the above definition reduces to the so-called  Fisher information
\begin{equation}\label{key}
I(\xi ) = {\int\limits_U {\left( {\frac{{\partial \ln f(\vec u; \xi )}}{{\partial \xi }}} \right)} ^2}f(\vec u;\xi ){\mathrm{d}^d}u.
\end{equation}

The Fisher metric plays the role of a Riemannian metric on the space of parameters $ \vec \xi  = ({\xi ^1},{\xi ^2}, \ldots ,{\xi ^n})$, where every point defines a different statistical model (or a PDF). We will not distinguish a given point $ \vec \xi$ in the parameter space and its associated PDF $f(\vec{u}; \vec \xi)$. Hence, given two points on the manifold their geodesic distance is interpreted as the statistical distinguishability of the PDFs \cite{costa:2015}.\\
\indent The action for the geodesics on the statistical manifold is given by the functional
\begin{equation}\label{eqGeodesicL}
L = \int\limits_{{r_1}}^{{r_2}} \sqrt {{g_{ab}}(\vec \xi )\,{\frac{\mathrm{d}\xi^a  (r)}{\mathrm{d}r}\frac{\mathrm{d}\xi^b(r)}{\mathrm{d}r}} } \,\mathrm{d}r,
\end{equation}
which under variation yields the system of geodesic equations
\begin{equation}\label{key1}
\frac{\mathrm{d}^2 \xi ^c(r)}{\mathrm{d}r^2} + \Gamma _{ab}^c \frac{\mathrm{d} \xi^a(r)}{\mathrm{d}r} \frac{\mathrm{d} \xi^b(r)}{\mathrm{d}r} = 0, \quad c=1,...,n.
\end{equation}
The invariant geodesic length $L$ between statistical models is then obtained from (\ref{eqGeodesicL}) after solving (\ref{key1}) for the geodesic profiles of the parameters $\xi^a(r)$ as functions of some proper ordering parameter $r$. 

For models with a single parameter one can determine the Fisher distance exactly up to a scale factor. For example, the Fisher information (metric) for the exponential distribution (\ref{eqExpD}) is given by
\begin{equation}\label{key}
g_{\lambda\lambda}=I(\lambda ) = \frac{1}{{{\lambda ^2}}}.
\end{equation}
Therefore, one can compute the distance function for this model directly by solving a single geodesic equation. For this purpose, we find the inverse metric  $g^{\lambda\lambda}=\lambda^2$ and the Christoffel symbol  $\Gamma^{\lambda}_{\lambda\lambda}=g^{\lambda\lambda}\partial_\lambda g_{\lambda\lambda}/2=-1/\lambda$. Thus, the geodesic equation for the model parameter $\lambda(r)$ is
\begin{equation}\label{key}
\frac{\mathrm{d}^2 \lambda}{\mathrm{d}r^2}  - \frac{1}{\lambda }\left(\frac{\mathrm{d} \lambda }{\mathrm{d}r}\right)^2 = 0
\end{equation}
with the simple solution $\lambda (r)=c_2 e^{c_1 r}$. Imposing boundary conditions, $\lambda(0)=\lambda_1$ and $\lambda (1)=\lambda_2$, one finds $c_1=\ln(\lambda_2/\lambda_1)$ and $c_2=\lambda_1$. Therefore, the geodesic distance between two statistical exponential models with corresponding parameters $\lambda_1$ and $\lambda_2$ is written by \cite{taylorE}
\begin{equation}\label{key}
{L_{12}} = L({\lambda _1},{\lambda _2}) = \left| {\int\limits_0^1 {\sqrt {{g_{\lambda \lambda }}\left(\frac{\mathrm{d}\lambda}{\mathrm{d}r}\right)^2} \mathrm{d}r} } \right| = \left| {\ln \frac{{{\lambda _2}}}{{{\lambda _1}}}} \right|.
\end{equation}

On the other hand, the Fisher metric for the log-normal distribution (\ref{eqLogND}) is given by
\begin{equation}\label{key}
d{s^2} = {g_{ab}}(\vec \xi )\mathrm{d}{\xi ^a}\mathrm{d}{\xi ^b} = \frac{{\mathrm{d}{\alpha ^2} + 2\mathrm{d}{\beta ^2}}}{{{\beta ^2}}}.
\end{equation}

\noindent The geodesic profiles for $\alpha(r)$ and $\beta(r)$ under this metric are given by the coupled system of ordinary second order differential equations
\begin{equation}\label{eqGeodLN}
\alpha ''(r) - \frac{{2\beta '(r)}}{{\beta (r)}}\alpha '(r) = 0,\quad \beta ''(r) - \frac{{\beta '{{(r)}^2}}}{{\beta (r)}} + \frac{{\alpha '{{(r)}^2}}}{{2\beta (r)}} = 0,
\end{equation}
together with the boundary conditions $\alpha(0)=\alpha_1$, $\alpha(1)=\alpha_2$, $\beta(0)=\beta_1$, $\beta(1)=\beta_2$.

In what follows, we will compute the Fisher distances between our models in the given ranges and find out how dissimilar they are from each other. For shortness of notation we will use the following indices: 1 for PSO, 2 for SA, and 3 for rSA. 

We begin by computing the Fisher distances between our exponential distributions for the time data in the range $R_3$, namely
\begin{equation}\label{key}
{L_{12}} = \left| {\ln \frac{{8.80}}{{14.30}}} \right| = 0.49,\quad {L_{13}} = \left| {\ln \frac{{10.01}}{{14.30}}} \right| = 0.36,\quad {L_{23}} = \left| {\ln \frac{{10.01}}{{8.00}}} \right| = 0.13.
\end{equation}
With similar computation one finds the Fisher distances in the range $R_4$:
\begin{equation}\label{key}
{L_{12}} = \left| {\ln \frac{{0.305}}{{0.162}}} \right| = 0.635,\quad {L_{13}} = \left| {\ln \frac{{0.305}}{{0.324}}} \right| = 0.060,\quad {L_{23}} = \left| {\ln \frac{{0.324}}{{0.162}}} \right| = 0.695.
\end{equation}
Finally, in the range $R_5$, one finds
\begin{equation}\label{key}
{L_{12}} = \left| {\ln \frac{{0.008}}{{0.002}}} \right| = 1.550,\quad {L_{13}} = \left| {\ln \frac{{0.008}}{{0.005}}} \right| = 0.512,\quad {L_{23}} = \left| {\ln \frac{{0.005}}{{0.002}}} \right| = 1.038.
\end{equation}
%

If we want to compare the log-normal models, we can find numerically the functions $\alpha(r)$ and $\beta(r)$ from (\ref{eqGeodLN}) and consequently  calculate the following integral:
\begin{equation}\label{key}
L_{ij} =\int\limits_0^1 {\frac{1}{{\beta(r)}}\sqrt {{{ \alpha }^{\prime2}} + 2{{\beta }^{\prime2}}} } \,\mathrm{d}r, \quad i,j=1,2,3, \quad i\ne j,
\end{equation}
with the proper boundary conditions, namely $\left( \alpha(0),\beta(0)\right)=\left(\alpha_i,\beta_i\right)$ and $\left( \alpha(1),\beta(1)\right)=\left(\alpha_j,\beta_j\right)$. In this case, one can show that in $R_3$ the Fisher distances between the tree log-normal models are
\begin{equation}\label{key}
{L_{12}} = 0.34,\quad {L_{13}} = 0.27,\quad {L_{23}} = 0.09.
\end{equation}
For the models in the mid range $R_4$ we find
\begin{equation}\label{key}
{L_{12}} = 0.40,\quad {L_{13}} = 0.20,\quad {L_{23}} = 0.46.
\end{equation}
And finally, in $R_5$, one has
\begin{equation}\label{key}
{L_{12}} = 0.88,\quad {L_{13}} = 0.67,\quad {L_{23}} = 0.81.
\end{equation}

It is useful to collect the results in tables (table \ref{distances}).
\begin{table}[h!]
	\small\centering
	\begin{subtable}{1.0\textwidth}
		\centering
		\begin{tabular}{||c |c |c| c||} 
			\hline
			\multicolumn{4} {||c||} {$R_3$} \\ 
			\hline
			& $L_{12}$ & $L_{13}$  & $L_{23}$ \\ 
			\hline
			Exponential & $0.49$ & $0.36$  & $0.13$ \\
			\hline
			Log-normal & $0.34$  & $0.27$ & $0.09$ \\
			\hline
		\end{tabular}
		\caption{$L_{12}$, $L_{13}$ and $L_{23}$ for the exponential and log-normal fits in the range $R_3$.\label{distancesR3}}
	\end{subtable}
	\vspace{1em}
	
	\begin{subtable}{1.0\textwidth}
		\centering
		\begin{tabular}{||c |c |c| c||} 
			\hline
			\multicolumn{4} {||c||} {$R_4$} \\ 
			\hline
			& $L_{12}$  & $L_{13}$& $L_{23}$ \\ 
			\hline
			Exponential & $0.635$ & $0.060$ & $0.695$ \\
			\hline
			Log-normal & $0.40$ & $0.20$ & $0.46$ \\
			\hline
		\end{tabular}
		\caption{$L_{12}$, $L_{13}$ and $L_{23}$ for the exponential and log-normal fits in the range $R_4$.\label{distancesR4}}
	\end{subtable}
	\vspace{1em}
	
	\begin{subtable}{1.0\textwidth}
		\centering
		\begin{tabular}{||c |c |c| c||} 
			\hline
			\multicolumn{4} {||c||} {$R_5$} \\ 
			\hline
			& $L_{12}$  & $L_{13}$ & $L_{23}$ \\ 
			\hline
			Exponential & $1.550$ & $0.512$ & $1.038$ \\
			\hline
			Log-normal & $0.88$ & $0.67$ & $0.81$ \\
			\hline
		\end{tabular}
		\caption{$L_{12}$, $L_{13}$ and $L_{23}$ for the exponential and log-normal fits in the range $R_5$.\label{distancesR5}}
	\end{subtable}
	\caption{Geodesic distances in the parameter spaces of the respective distributions between the three algorithms in the tested ranges. We use the indices of $L$ to denote the following: 1 for PSO, 2 for SA, and 3 for rSA. \label{distances}}
\end{table}

One can infer that in the lowest range $R_3$, when considering the exponential distribution, the SA and rSA algorithms are similar relative to each other ($L_{23}=0.13$, i.e. they are closest), while they are quite dissimilar to PSO ($L_{12}=0.49$ and $L_{13}=0.36$). The same is valid also for the log-normal model in $R_3$.

On the other hand, in the mid range $R_4$, the PSO and rSA algorithms are similar relative to each other, for example $L_{13}=0.060$, while they are notably dissimilar to SA, i.e. $L_{12}=0.635$ and $L_{23}=0.695$. This result persists also in the next range $R_5$.
 
\section{Conclusion}\label{concl}

\indent In this paper we adapted the number-theoretic sum of three cubes problem to an optimisation setting. This was motivated by the desire to use some random search algorithm to hopefully improve the time it takes to find a solution. Turning the problem into an optimisation one was not hard and resulted in equation (\ref{fdef}). However, finding a global minimum to (\ref{fdef}) with sufficient speed turned out to be a highly non-trivial task (as was expected).\\
\indent Our attempts in this direction led us to test the performace of three different search heuristics in three ranges for $(x,y)$ when applied to our problem in the special case $k=2$. The first one is loosely based on DFO with some major modifications, while the second and third one are more or less direct implementations of SA and rSA, respectively.\\
\indent The metric for the performance of the algorithms was the time it takes a given algorithm to reach a solution to (\ref{deq}) (i.e. a global minimum of (\ref{fdef})). After a large number of tests, we analysed the results by fitting the respective datasets of running times with two different distributions -- exponential and log-normal.\\
\indent We have analysed two specific aspects of the algorithms, namely their time performance and their similarity. A conclusion about the time performance can be made by looking both at the mean and the mode of the running times (collected in table \ref{statchars}), while the relative similarity between the algorithms can be measured by the Fisher distances between the respective PDFs (table \ref{distances}).\\
\begin{table}[ht]
	\small\centering
	\begin{subtable}{1.0\textwidth}
		\centering
		\begin{tabular}{||c |c |c| c||} 
			\hline
			\multicolumn{4} {||c||} {$R_3$} \\ 
			\hline
			& PSO & SA & rSA \\ 
			\hline
			$\bar{t}_{\text{exp}}\ [s]$ & $0.070\pm 0.001$ & $0.114\pm 0.002$ & $0.100\pm 0.002$ \\
			\hline
			$\bar{t}_{\text{l-n}} \ [s]$ & $0.089\pm 0.003$ & $0.151\pm 0.006$ & $0.130\pm 0.005$ \\
			\hline
			$\mathrm{Mode}_{\text{l-n}}[t] \ [s]$ & $0.008$ & $0.010$ & $0.010$ \\ 
			\hline
		\end{tabular}
		\caption{$\bar{t}_{\text{exp}}$, $\bar{t}_{\text{l-n}}$ and $\mathrm{Mode}_{\text{l-n}}[t]$ for PSO, SA and rSA in the range $R_3$.\label{statcharsR3}}
	\end{subtable}
\end{table}
\vspace{-1.2em}
	\begin{table}[ht]\ContinuedFloat
	\small\centering
	\begin{subtable}{1.0\textwidth}
		\centering
		\begin{tabular}{||c |c |c| c||} 
			\hline
			\multicolumn{4} {||c||} {$R_4$} \\ 
			\hline
			& PSO & SA & rSA \\ 
			\hline
			$\bar{t}_{\text{exp}}\ [s]$ & $3.27\pm 0.07$ & $6.18\pm 0.12$ & $3.08\pm 0.06$ \\
			\hline
			$\bar{t}_{\text{l-n}} \ [s]$ & $4.17\pm 0.14$ & $8.86\pm 0.38$ & $4.14\pm 0.17$ \\
			\hline
			$\mathrm{Mode}_{\text{l-n}}[t] \ [s]$ & $0.36$ & $0.32$ & $0.21$ \\ 
			\hline
		\end{tabular}
		\caption{$\bar{t}_{\text{exp}}$, $\bar{t}_{\text{l-n}}$ and $\mathrm{Mode}_{\text{l-n}}[t]$ for PSO, SA and rSA in the range $R_4$.\label{statcharsR4}}
	\end{subtable}
	\end{table}
	\vspace{-1.2em}
\begin{table}[ht]\ContinuedFloat
\small\centering
	\begin{subtable}{1.0\textwidth}
		\centering
		\begin{tabular}{||c |c |c| c||} 
			\hline
			\multicolumn{4} {||c||} {$R_5$} \\ 
			\hline
			& PSO & SA & rSA \\ 
			\hline
			$\bar{t}_{\text{exp}}\ [s]$ & $123.8\pm 8.2$  & $583.5\pm 38.6$ & $206.6\pm 13.7$ \\
			\hline
			$\bar{t}_{\text{l-n}} \ [s]$ & $154.5\pm 17.0$ & $847.9\pm 141.4$ & $421.3\pm 101.3$ \\
			\hline
			$\mathrm{Mode}_{\text{l-n}}[t] \ [s]$ & $15.1$ & $16.5$ & $1.1$ \\ 
			\hline
		\end{tabular}
		\caption{$\bar{t}_{\text{exp}}$, $\bar{t}_{\text{l-n}}$ and $\mathrm{Mode}_{\text{l-n}}[t]$ for PSO, SA and rSA in the range $R_5$. \label{statcharsR5}}
	\end{subtable}
	\caption{Expected values and modes for the respective distribution fits for all the algorithms in the tested ranges.\label{statchars}}
\end{table}

\indent The main conclusion, when considering the average times, is that for this particular problem our version of PSO performs similarly to rSA in the range $R_4$, but better in the ranges $R_3$ and $R_5$. As expected rSA is better than SA in all ranges. The relative performance of the algorithms in the different ranges as measured by the ratios of the average times is shown in table \ref{relPerform}.\\
\begin{table}[h!]
	\small\centering
\begin{tabular}{ ||c| c| c| c| c |c|c|c|c||} 
\hline
Exponential &$R_3$ &$R_4$  &$R_5$& &Log-normal &$R_3$ &$R_4$  &$R_5$  \\	\hline
$\bar t_{SA}/\bar t_{PSO}$ &1.6 &1.9 &4.7& &$\bar t_{SA}/\bar t_{PSO}$ &1.7 &2.1 &5.5\\	\hline
$\bar t_{rSA}/\bar t_{PSO}$ &1.4 &0.9 &1.7& &$\bar t_{rSA}/\bar t_{PSO}$ &1.5 &1.0 &2.7\\	\hline
$\bar t_{SA}/\bar t_{rSA}$ &1.1 &2.0 &2.8& &$\bar t_{SA}/\bar t_{rSA}$ &1.1 &2.1 &2.0\\	\hline
\end{tabular}
	\caption{Relative performance of the algorithms in the tested ranges. For example, the exponential model fit  in $R_3$ shows that our PSO is 1.6 times faster on average than SA and 1.4 times faster than rSA. On the other hand, the log-normal fit states that PSO is 1.7 times faster than SA and 1.5 times faster than rSA on average. \label{relPerform}}
\end{table}

\indent When looking at the mode of the respective log-normal distributions, we see a slightly different picture in the highest range -- rSA is by far the best method, which is not so pronounced in the lower ranges. This is evident from the distribution of its running times, which has $\mathrm{Mode}[t]=1.1\,s$ in $R_5$, compared to $15.1\,s$ and $16.5\,s$ for PSO and SA respectively (table \ref{statchars}). In other words, rSA finds most of the solutions notably more quickly than PSO and SA.\\
\indent Finally, considering the Fisher distances between the respective distributions for the different algorithms, we can see that in the lowest range SA and rSA are the most similar. This changes in the higher ranges, where the distance between PSO and rSA is the smallest.

\indent Now, where can we go from here? The final goal is to be able to find solutions to (\ref{deq}) for various values of $k$ in ranges above $10^{20}$ in reasonable time. As is well known, the solution density there is significantly reduced and this means that the search becomes very time consuming. We believe that some stochastic algorithm can be found that can produce solutions in acceptable time.\\
\indent One way to reduce the search time is parallelisation. Generally, this can be done in many ways. As was already mentioned, the mode of the running times of rSA in the highest range is peculiarly small. This suggests probably the simplest method to achieve some sort of parallelisation -- run the same instance of the algorithm on many cores and just wait for the first one to finish. The probability of achieving a running time in $R_5$ with rSA of less than $1.5s$ for example is $P_{\text{rSA}}\left(t\le 1.5\right)\approx 0.0333$. Compare that to the probabilities for the same event with the other two algorithms: $P_{\text{SA}}\left(t\le 1.5\right)\approx 0.00096$ and $P_{\text{PSO}}\left(t\le 1.5\right)\approx 0.00097$.\\
\indent Another line of investigation is to search for a better heuristic. In general, stochastic optimisation algorithms are highly specific to the problem and finding a good one isn't easy. A promising new development with regards to this is \cite{li16}, which may enable us to delegate the task to AI.

\subsection*{Acknowledgements}

The authors would like to thank R. C. Rashkov, S. Yazadjiev, H. Dimov, P. Nedkova, G. Gyulchev and K. Staykov for many useful comments and discussions. T. V. gratefully acknowledges the support of the Bulgarian national program
“Young Scientists and Postdoctoral Research Fellows 2020”, and the Sofia University Grant 80-10-68. This work was also partially supported by the Bulgarian NSF grant N28/5.

\newpage
\begin{appendices}
	
	\section{ C code implementation}
	
	\subsection{PSO}\label{psocode}
	{\tiny 
\begin{lstlisting}[language=C,breaklines=true]
#include <stdio.h>
#include <stdlib.h>
#include <math.h>
#include <time.h>

double ff(int x, int y, int kk); //The fitness function.
int sort_and_break(double q[], int n); //A function which chooses the number of a best particle breaking ties randomly.



int main()
{
//Parameters:
        int kk = 2; //Right side of the diophantine equation.
        int lb[2] = {-pow(10,4),0}; //Lower bounds for the search space.
        int ub[2] = {0,pow(10,4)}; //Upper bounds for the search space.
        double thr = pow(10,-5); //Solution threshold. Used to decide when to check whether a solution has been found.
        int s = 50; //Number of particles.
        int dp = (s/5); //Number of particles that need to be at the best swarm best position in order to disperse all particles.

//Variables:
	time_t t0 = time(NULL); //Time used to seed the RNG.
	int pos[s][2]; //Particle positions.
	double fitf[s]; //Fitness functions for the respective particles.
	int nbrs[s][3]; //Neighbours of each particle.
	int sb[2],bsb[2],nb[s][2]; //Position of the swarm best particle, best swarm best position and position of the neighbours best particles.
	double sbfitf,bsbfitf,nbfitf[s]; //Fitness functions for the above.
        int dc; //Dispersion counter to decide when to disperse the particles.
	double nbrsfitf[3],r; //Auxiliary double variables.
        int i,j,zt; //Auxiliary integer variables.

//Seeding the RNG:
	srand48(t0);

//Calculating the neighbours of each particle:
	i = 0; while(i < s)
	{
		nbrs[i][0] = (s+i-1)%s;
		nbrs[i][1] = (s+i)%s;
		nbrs[i][2] = (s+i+1)%s;
		i = i+1;
	}

//Initialisation:
	i = 0; while(i < s)
	{
		pos[i][0] = lround(drand48()*(ub[0]-lb[0])+lb[0]);
		pos[i][1] = lround(drand48()*(ub[1]-lb[1])+lb[1]);
		fitf[i] = ff(pos[i][0],pos[i][1],kk);
		i = i+1;
	}

//Determining the swarm best and best swarm best particles:
	i = sort_and_break(fitf, s);
	sb[0] = pos[i][0]; sb[1] = pos[i][1]; sbfitf = fitf[i];
	bsb[0] = sb[0]; bsb[1] = sb[1]; bsbfitf = sbfitf;

	if(bsbfitf <= thr) //Solution check.
	{
		zt = lround(cbrt(kk-pow(bsb[0],3)-pow(bsb[1],3)));
		if(pow(bsb[0],3)+pow(bsb[1],3)+pow(zt,3) == kk)
		{
			printf("{%d,%d,%d}\n",bsb[0],bsb[1],zt);
			goto end;
		}
	}

//Calculating the best particle among the neighbours:
	i = 0; while(i < s)
	{
		nbrsfitf[0] = fitf[nbrs[i][0]];
		nbrsfitf[1] = fitf[nbrs[i][1]];
		nbrsfitf[2] = fitf[nbrs[i][2]];
		j = sort_and_break(nbrsfitf, 3);
		nb[i][0] = pos[(s+i+j-1)%s][0];
		nb[i][1] = pos[(s+i+j-1)%s][1];
		nbfitf[i] = fitf[(s+i+j-1)%s];
		i = i+1;
	}

//Iterations:
	while(1)
	{
//Position update:
		if(dc >= dp)
		{
			i = 0; while(i < s)
			{
				j = 0; while(j < 2)
				{
					r = drand48();
					pos[i][j] = lround(r*(ub[j]-lb[j])+lb[j]);
					j = j+1;
				}
				i = i+1;
			}
		}
		else
		{
			i = 0; while(i < s)
			{
				j = 0; while(j < 2)
				{
					r = drand48();
					pos[i][j] = lround(nb[i][j]+(r/2.0)*(sb[j]+bsb[j]-2.0*pos[i][j]));
					j = j+1;
				}
				i = i+1;
			}
		}

//Confinement and dispersion condition check:
		i = 0; dc = 0; while(i < s)
		{
			if((pos[i][0] < lb[0]) || (pos[i][0] > ub[0]) || (pos[i][1] < lb[1]) || (pos[i][1] > ub[1]))
			{
				pos[i][0] = lround(drand48()*(ub[0]-lb[0])+lb[0]);
				pos[i][1] = lround(drand48()*(ub[1]-lb[1])+lb[1]);
			}
			fitf[i] = ff(pos[i][0],pos[i][1],kk);
			if((pos[i][0] == bsb[0]) && (pos[i][1] == bsb[1])) { dc = dc+1; }
			i = i+1;
		}

//Swarm best, best swarm best and neighbours best updates:
		i = sort_and_break(fitf, s);
		sb[0] = pos[i][0]; sb[1] = pos[i][1]; sbfitf = fitf[i];
		r = drand48();
		if(sbfitf < bsbfitf)
		{
			bsb[0] = sb[0]; bsb[1] = sb[1]; bsbfitf = sbfitf;

			if(bsbfitf <= thr) //Solution check.
			{
				zt = lround(cbrt(kk-pow(bsb[0],3)-pow(bsb[1],3)));
				if(pow(bsb[0],3)+pow(bsb[1],3)+pow(zt,3) == kk)
				{
					printf("{%d,%d,%d}\n",bsb[0],bsb[1],zt);
					goto end;
				}
			}

		}

		if((sbfitf > bsbfitf) && (r <= (1-((sbfitf-bsbfitf)/0.5))))
		{
			bsb[0] = sb[0]; bsb[1] = sb[1]; bsbfitf = sbfitf;

			if(bsbfitf <= thr) //Solution check.
			{
				zt = lround(cbrt(kk-pow(bsb[0],3)-pow(bsb[1],3)));
				if(pow(bsb[0],3)+pow(bsb[1],3)+pow(zt,3) == kk)
				{
					printf("{%d,%d,%d}\n",bsb[0],bsb[1],zt);
					goto end;
				}
			}

		}

//Calculating the best particle among the neighbours:
		i = 0; while(i < s)
		{
			nbrsfitf[0] = fitf[nbrs[i][0]];
			nbrsfitf[1] = fitf[nbrs[i][1]];
			nbrsfitf[2] = fitf[nbrs[i][2]];
			j = sort_and_break(nbrsfitf, 3);
			nb[i][0] = pos[(s+i+j-1)%s][0];
			nb[i][1] = pos[(s+i+j-1)%s][1];
			nbfitf[i] = fitf[(s+i+j-1)%s];
			i = i+1;
		}
	}

	end:

	return 0;
}



double ff(int x, int y, int kk)
{
	double z;

	z = fabs(cbrt(kk-pow(x,3)-pow(y,3))-lround(cbrt(kk-pow(x,3)-pow(y,3))));

	return z;
}

int sort_and_break(double q[], int n)
{
        int k,i;
        int c = 0;
        int b[n];
        double gv = q[0];

	i = 0; while(i < n) { b[i] = 0; i = i+1; }

        i = 1; while(i < n)
        {
                if(q[i] == gv)
                {
                        c = c+1;
                        b[c] = i;
                }
                if(q[i] < gv)
                {
                        gv = q[i];
                        c = 0;
                        b[c] = i;
                }
                i = i+1;
        }

        k = lround(drand48()*c);

        return b[k];
}
\end{lstlisting}}

	\subsection{SA}\label{sacode}
	{\tiny 
\begin{lstlisting}[language=C,breaklines=true]
#include <stdio.h>
#include <math.h>
#include <stdlib.h>
#include <time.h>

double ff(int x,int y,int kk); //The energy function.
double temp(int m); //The cooling schedule (temperature).



int main()
{
//Parameters:
        int kk = 2; //Right side of the diophantine equation.
        int lbx = -pow(10,4); //Lower bound for the x coordinate of a state.
        int ubx = 0; //Upper bound for the x coordinate.
        int lby = 0; //Lower bound for the y coordinate of a state.
        int uby = pow(10,4); //Upper bound for the y coordinate.
        double thr = pow(10,-5); //Solution threshold. Used to decide when to check whether a solution has been found.

//Variables:
	time_t t0 = time(NULL); //Time used to seed the RNG.
	int x,y,xn,yn; //Coordinates of the current and new states.
	double z,zn; //Energy values for the above.
	int m; //Time.
	int a,b,c,d; //Variables used in the generation of a new state.
	double prob; //Probability for accepting a transition to the already generated state.
	int zt; //Auxiliary integer varialbe.
	double r; //Auxiliary double variable.

//Seeding the RNG:
	srand48(t0);

//Random initial state and its energy:
	x = lround(drand48()*(ubx-lbx)+lbx); y = lround(drand48()*(uby-lby)+lby);
	z = ff(x,y,kk);

	if(z <= thr) //Solution check.
	{
		zt = lround(cbrt(kk-pow(x,3)-pow(y,3)));
		if(pow(x,3)+pow(y,3)+pow(zt,3) == kk)
		{
			printf("{%d,%d,%d}\n",x,y,zt);
			goto end;
		}
	}

//Iterations:
	m = 1; while(1)
	{
//Confinement:
		a = -10; b = 10; c = -10; d = 10;
		if(x+a <= lbx) { a = lbx-x; }
		if(x+b >= ubx) { b = ubx-x; }
		if(y+c <= lby) { c = lby-y; }
		if(y+d >= uby) { d = uby-y; }

//State generation:
		roll:
		xn = x+lround(drand48()*(b-a)+a);
		yn = y+lround(drand48()*(d-c)+c);
		if((xn == x) && (yn == y)) { goto roll; }
		zn = ff(xn,yn,kk);

//State transition:
		if(zn <= z)
		{
			x = xn;
			y = yn;
			z = zn;
		}
		else
		{
			prob = exp((z-zn)/temp(m));
			r = drand48();
			if(r <= prob)
			{
				x = xn;
				y = yn;
				z = zn;
			}
		}

		if(z <= thr) //Solution check.
		{
			zt = lround(cbrt(kk-pow(x,3)-pow(y,3)));
			if(pow(x,3)+pow(y,3)+pow(zt,3) == kk)
			{
				printf("{%d,%d,%d}\n",x,y,zt);
				goto end;
			}
		}
		m = m+1;
	}

	end:

	return 0;
}



double ff(int x,int y,int kk)
{
	double z;

	z = fabs(cbrt(kk-pow(x,3)-pow(y,3))-lround(cbrt(kk-pow(x,3)-pow(y,3))));

	return z;
}



double temp(int m)
{
	double z;

	z = 1.0/(log(m)+0.01);

	return z;
}
\end{lstlisting}}

	\subsection{rSA}\label{rsacode}
	{\tiny 
\begin{lstlisting}[language=C,breaklines=true]
#include <stdio.h>
#include <math.h>
#include <stdlib.h>
#include <time.h>

double ff(int x,int y,int kk); //The energy function.
double temp(int m); //The cooling schedule (temperature).



int main()
{
//Parameters:
	int kk = 2; //Right side of the diophantine equation.
	int lbx = -pow(10,4); //Lower bound for the x coordinate of a state.
	int ubx = 0; //Upper bound for the x coordinate.
	int lby = 0; //Lower bound for the y coordinate of a state.
	int uby = pow(10,4); //Upper bound for the y coordinate.
	double thr = pow(10,-5); //Solution threshold. Used to decide when to check whether a solution has been found.
        int rtm = 30; //Number of consecutive states with equal energies needed for a restart.

//Variables:
        time_t t0 = time(NULL); //Time used to seed the RNG.
        int xo,yo,x,y,xn,yn; //Coordinates of the old, current and new states.
	double zo,z,zn; //Energy values for the above.
	int m; //Time.
	int a,b,c,d; //Variables used in the generation of a new state.
	double prob; //Probability for accepting a transition to the already generated state.
	int rt; //Current number of consecutive states with equal energies.
	int zt; //Auxiliary integer variable.
	double r; //Auxiliary double variable.

//Seeding the RNG:
	srand48(t0);

	restart:

//Random initial state and its energy:
	rt = 1;
	x = lround(drand48()*(ubx-lbx)+lbx); y = lround(drand48()*(uby-lby)+lby);
	z = ff(x,y,kk);

	if(z <= thr) //Solution check.
	{
		zt = lround(cbrt(kk-pow(x,3)-pow(y,3)));
		if(pow(x,3)+pow(y,3)+pow(zt,3) == kk)
		{
			printf("{%d,%d,%d}\n",x,y,zt);
			goto end;
		}
	}

//Iterations:
	m = 1; while(1)
	{
//Confinement:
		a = -10; b = 10; c = -10; d = 10;
		if(x+a <= lbx) { a = lbx-x; }
		if(x+b >= ubx) { b = ubx-x; }
		if(y+c <= lby) { c = lby-y; }
		if(y+d >= uby) { d = uby-y; }

//State generation:
		roll:
		xn = x+lround(drand48()*(b-a)+a);
		yn = y+lround(drand48()*(d-c)+c);
		if((xn == x) && (yn == y)) { goto roll; }
		zn = ff(xn,yn,kk);

//State transition:
		xo = x;
		yo = y;
		zo = z;
		if(zn <= z)
		{
			x = xn;
			y = yn;
			z = zn;
		}
		else
		{
			prob = exp((z-zn)/temp(m));
			r = drand48();
			if(r <= prob)
			{
				x = xn;
				y = yn;
				z = zn;
			}
		}

//Restart condition check:
		if(zo == z)
		{
			rt = rt+1;
			if(rt == rtm) { goto restart; }
		}
		else
		{
			rt = 1;
		}

		if(z <= thr) //Solution check.
		{
			zt = lround(cbrt(kk-pow(x,3)-pow(y,3)));
			if(pow(x,3)+pow(y,3)+pow(zt,3) == kk)
			{
				printf("{%d,%d,%d}\n",x,y,zt);
				goto end;
			}
		}
		m = m+1;
	}

	end:

	return 0;
}



double ff(int x,int y,int kk)
{
	double z;

	z = fabs(cbrt(kk-pow(x,3)-pow(y,3))-lround(cbrt(kk-pow(x,3)-pow(y,3))));

	return z;
}



double temp(int m)
{
	double z;

	z = 1.0/(log(m)+0.01);

	return z;
}
\end{lstlisting}}

\end{appendices}

\end{document}